\documentclass[12pt]{article}
\usepackage{amssymb}
\usepackage{a4}
\begin{document}
\def\st{\, : \,}
\def\kbar{{\mathchar'26\mkern-9muk}}  
\def\bra#1{\langle #1 \vert}
\def\ket#1{\vert #1 \rangle}
\def\vev#1{\langle #1 \rangle}
\def\ad{\mbox{ad}\,}
\def\ker{\mbox{Ker}\,}
\def\im{\mbox{Im}\,}
\def\der{\mbox{Der}\,}
\def\ad{\mbox{ad}\,}
\def\b#1{{\mathbb #1}}
\def\c#1{{\cal #1}}
\def\pt{\partial_t}
\def\px{\partial_1}
\def\bpx{\bar\partial_1}
\def\la{\langle}
\def\ra{\rangle}
\def\nn{\nonumber \\}
\def\pn{\par\noindent}
\def\etal{{\it et al.}\ }
\def\sq{\mbox{\rlap{$\sqcap$}$\sqcup$}}
\def\R{{\cal R}\,}
\newcommand{\tr}{\triangleright\,}
\newcommand{\tl}{\,\triangleleft}
\newcommand{\tro}{\triangleright^{op}\,}
\newcommand{\tlo}{\,\stackrel{op}{\triangleleft}}
\def\cross{{\triangleright\!\!\!<}}
\def\cocross{{>\!\!\!\triangleleft\,}}
\def\uqg{\mbox{$U_q{\/\bf g}$ }}
\def\uqso{\mbox{$U^{op}_q{\/so(N)}$ }}
\def\uqsp{\mbox{$U_q^+{\/so(N)}$ }}
\def\uqsn{\mbox{$U_q^-{\/so(N)}$ }}
\def\uqs{\mbox{$U_q{\/so(N)}$ }}

\renewcommand{\thefootnote}{\fnsymbol{footnote}}

\renewcommand{\theequation}{\thesection.\arabic{equation}}
\newcommand{\initiate}{\setcounter{equation}{0}}

\newcommand{\ba}{\begin{array}}
\newcommand{\ea}{\end{array}}
\newcommand{\be}{\begin{equation}}
\newcommand{\ee}{\end{equation}}
\newcommand{\bea}{\begin{eqnarray}}
\newcommand{\eea}{\end{eqnarray}}
\newcommand{\beas}{\begin{eqnarray*}}
\newcommand{\eeas}{\end{eqnarray*}}
%
%
%
\newtheorem{prop}{Proposition}
\newtheorem{lemma}{Lemma}
\newtheorem{theorem}{Theorem}
\newtheorem{corollary}{Corollary}
%
%
\newenvironment{proof}[1]{\vspace{5pt}\noindent{\bf Proof #1}\hspace{6pt}}%
{\hfill\sq}
\newcommand{\bp}{\begin{proof}}
\newcommand{\ep}{\end{proof}\par\vspace{10pt}\noindent}
%
%

\title{Geometrical Tools for Quantum Euclidean Spaces}

\author{B.L. Cerchiai,$\strut^{1,2}$ \,
        G. Fiore,$\strut^{3,4}$ \, J. Madore$\strut^{5,2}$ \\\\
        \and
        $\strut^1$Sektion Physik, Ludwig-Maximilian Universit\"at,\\
        Theresienstra\ss e 37, D-80333 M\"unchen
        \and
        $\strut^2$Max-Planck-Institut f\"ur Physik\\
        F\"ohringer Ring 6, D-80805 M\"unchen
        \and
        $\strut^3$Dip. di Matematica e Applicazioni, Fac.  di Ingegneria\\ 
        Universit\`a di Napoli, V. Claudio 21, 80125 Napoli
        \and
        $\strut^4$I.N.F.N., Sezione di Napoli,\\
        Mostra d'Oltremare, Pad. 19, 80125 Napoli
        \and
        $\strut^5$Laboratoire de Physique Th\'eorique et Hautes Energies\\
        Universit\'e de Paris-Sud, B\^atiment 211, F-91405 Orsay
        }
\date{}

\maketitle
\abstract{We apply one of the formalisms of noncommutative geometry to
$\b{R}^N_q$, the quantum space covariant under the quantum group
$SO_q(N)$. Over $\b{R}^N_q$ there are two $SO_q(N)$-covariant
differential calculi. For each we find a frame, a metric and two
torsion-free covariant derivatives which are metric compatible up to a
conformal factor and which have a vanishing linear curvature.
This generalizes results found in a previous article
for the case of $\b{R}^3_q$. As in the case $N=3$, one has to 
slightly enlarge the algebra $\b{R}^N_q$;
for $N$ odd one needs only one new generator whereas for $N$ 
even one needs two. As in the particular case $N=3$ there is 
a conformal ambiguity in the natural metrics on the 
differential calculi over $\b{R}^N_q$. While in our previous 
article the frame was found `by hand', here we disclose the 
crucial role of the quantum group covariance and exploit it in 
the construction. As an intermediate step, we find a 
homomorphism from the cross product of $\b{R}^N_q$ with $U_qso(N)$ 
into $\b{R}^N_q$, an interesting result in itself.
}

\parskip 4pt plus2pt minus2pt
\noindent
LMU-TPW 99-17, MPI-PhT/99-45\\
\noindent
Preprint 99-52 Dip. Matematica e Applicazioni, Universit\`a di Napoli
\newpage

\initiate
\section{Introduction}

It is an old idea \cite{Sny47a, Sny47b} that a noncommutative
modification of the algebraic structure of space-time could provide a
regularization of the divergences of quantum field theory, because the
representations of noncommutative `spaces' have a lattice-like
structure~\cite{HebSchSchWeiWes94, Fio95, KemMan97, FicLorWes96,
CerHinMadWes99a} which should automatically impose an ultraviolet
cut-off~\cite{DopFreRob95, ChoHinMadSte99}. This idea has been
challenged recently from various points of view~\cite{Fil96,
ChaDemPre98, MadSte99c, MarSan99, SeiWit99}.  In any case to discuss
it and other problems it is necessary to have a noncommutative
version of flat space~\cite{FadResTak89, OgiSchWesZum92, LukNowRue92,
Cas95, CerWes98}.  Two approaches have been suggested to endow an
algebra with the noncommutative generalization of a differential
calculus, that of Woronowicz~\cite{Wor79, Wor87a, Wor87b} and that of
Connes~\cite{Con86, Con94}. The formalism which we shall use here is
an attempt to conciliate these two points of view in a particular
class of examples.  To this end we use a particular noncommutative
version \cite{DimMad96} of the moving-frame formalism of E. Cartan.

In a previous article a detailed study was made of the noncommutative
geometry of $\b{R}^3_q$, the quantum space covariant under the quantum
group $SO_q(3)$. It was found that one had to slightly extend the
algebra by adding a `dilatator' $\Lambda$ in order to reduce the
center to the complex numbers and so to be able to construct an
essentially unique metric, whereas for the construction of a frame
one also had to add the square roots and the inverses of some
generators.

The results are here extended to the case of the general algebra
$\b{R}^N_q$. We find that the cases $N$ even and $N$ odd are somewhat
different. When $N$ is odd the formalism is quite similar to the case
$N=3$. When, on the other hand, $N$ is even, yet another extension
must be made for the construction of a frame.  We must add to the
algebra one of the components $K$ of the angular momentum in order to
have a trivial center. The differential calculus of
\cite{CarSchWat91b} is extended by setting $dK = 0$ and either
$d\Lambda = 0$ or $d\Lambda \neq 0$ but fixed by a modified Leibniz
rule.  These extensions are in a sense unsatisfactory since they imply
that there are elements of the extended algebra (what we shall call
$\c{A}_N$) which have vanishing derivative but which are none-the-less
noncommutative analogues of non-constant functions.  Their inclusion
can be interpreted as an embedding of the `configuration space' into
part of `phase space'. The different possibilities lead to a conformal
ambiguity in the natural metrics on the differential calculi over
$\b{R}^N_q$; one choice favours the geometry $S^{N-1}\times \b{R}$ and
the second, the one we emphasize here, favours the flat geometry
$\b{R}^N$.  For each of its two $SO_q(N)$-covariant differential
calculi we find the corresponding frame and two torsion-free covariant
derivatives that are metric compatible up to a conformal factor and
which yield both a vanishing linear curvature. Apart from a few notes,
we leave the study of the reality structures, $*$-representations and
of the commutative limit as subjects to be treated elsewhere.

In Section~2 we briefly recall the tools of noncommutative
geometry~\cite{Con94} which will be needed. We start with a formal
noncommutative algebra $\c{A}$ and with a differential calculus
$\Omega^*(\c{A})$ over it. We define then a frame or
`Stehbein'~\cite{DimMad96} and the corresponding metric and covariant
derivative.  We also recall how a generalized Dirac
operator~\cite{Con94} can be constructed from the frame and a dual set
of inner derivations. Finally we recall the compatibility condition
between the metric and the covariant derivative.  The frame is what
will permit us to pass from the covariant definition of a differential
calculus~\cite{WesZum90}, with its emphasis on $q$-deformed
commutators, to the definition of Connes, which uses ordinary
commutators, at least on a formal level; we do not attempt to discuss
the Dirac operator as an operator on a graded Hilbert space.

In Section~3 we recall some formulae from the pioneering work of
Faddeev, Reshetikhin and Takhtajan~\cite{FadResTak89} on the
definition of $\b{R}_q^N$ as given by the coaction of the quantum
group $SO_q(N)$.  We give then a brief overview of the work of
Carow-Watamura, Schlieker, Watamura~\cite{CarSchWat91b} and
Ogievetski~\cite{Ogi92a} on the construction of two differential
calculi on $\b{R}_q^N$, which are based on the $\hat R$-matrix
formalism and are covariant with respect to $SO_q(N)$. They both yield
the de Rham calculus in the commutative limit.  However here it is
convenient to formulate quantum group covariance in terms of the
action on $\b{R}_q^N$ of the dual Hopf algebra $U_qso(N)$, rather than
in terms of the coaction of $SO_q(N)$.

In Section~4 we proceed with the actual construction of the frame over
$\b{R}^N_q$ and of the inner derivations dual to it.  We first solve
the problem in a larger algebra, $\Omega^*(\c{A}_N) \cocross\uqs$,
where we show that a frame has to transform under the action of the
quantum group $\uqso$ with opposite coalgebra.  We also find a dual
set of inner derivations by decomposing in the frame basis the formal
`Dirac operator', which had already been found~\cite{Zum97, Ste96} 
previously. It would be interesting, but requires some nontrivial
handwork, to  compare our results with the ones of Ref. 
\cite{AscCast96,AscCastSca97}. There multiparametric deformations of
the inhomogeneous $SO_q(N)$ quantum groups are considered, whereby  
multiparametric deformations (including
as a particular case the one-parameter one at hand) of the  
Euclidean space are obtained
by projection. The frame of the quantum group in the Woronowicz 
bicovariant differential calculi sense, i.e. the left- (or right-) 
invariant 1-forms, might also be projected and compared to ours.

Then in Section~5 we show that it is possible to find homomorphisms
$\varphi^\pm:\c{A}_N\cocross U_q^{\pm}{\/so(N)} \rightarrow\c{A}_N$
which act trivially on the factor $\c{A}_N$ on the left-hand side and
which project the components of the frame and of the inner derivations
from elements of $\c{A}_N \cocross\uqs$ onto elements in $\c{A}_N$. This
implies that in the $x^i$ basis they satisfy the `RLL' and the `gLL'
relations fulfilled by the $L^{\pm}$~\cite{FadResTak89} generators of
$\uqs$.  In the case that $N$ is odd it is possible to `glue' the
homomorphisms together to an isomorphism from the whole of
$\c{A}_N\cocross \uqs$ to $\c{A}_N$, an interesting and surprising
result in itself.

Finally in Section~6 we see that for each of the two calculi there is
essentially a unique metric, and two torsion-free $SO_q(N)$-covariant
linear connections which are compatible with it up to a conformal
factor.

\initiate
\section{The Cartan formalism}

In this section we briefly review a noncommutative extension 
\cite{DimMad96} of the moving-frame formalism of E. Cartan. We start
with a formal noncommutative associative algebra $\c{A}$ with a differential
calculus $\Omega^*(\c{A})$. If $\c{A}$ has a commutative limit and if
this limit is the algebra of functions on a manifold $M$ then we
suppose that the limit of the differential calculus is the ordinary
$\Omega^*(\c{A}_N)\cocross\uqs$de~Rham differential calculus on $M$.
We shall concentrate on the
case where the module of the 1-forms $\Omega^1(\c{A})$ is free of rank
$N$ as a left or right module and admits a special basis
$\{ \theta^a \}_{1 \le a \le N}$, referred to as `frame' or
`Stehbein', which commutes with the elements of $\c{A}$:
\be
[f,\theta^a]=0.                                              \label{thetaf}
\ee
This means that if the limit manifold exists it must be parallelizable.
The integer $N$ plays the role of the dimension of the manifold.
We suppose further that the basis $\theta^a$ is dual to a set of inner
derivations $e_a=\ad \lambda_a$ such that:
\be
df=e_a f \theta^a=[\lambda_a,f] \theta^a                         \label{df}
\ee
for any $f \in \c{A}$. The formal `Dirac operator'~\cite{Con94}, 
defined by the equation
\be
df = -[\theta,f],                                              \label{dirac}
\ee
is then given by
\be
\theta = -\lambda_a \theta^a.
\ee
We shall consider only the case where the center $\c{Z}(\c{A})$ of
$\c{A}$ is trivial: $\c{Z}(\c{A})=\b{C}$. If the original algebra does
not have a trivial center then we shall extend it an algebra which
does. The (wedge) product $\pi$ in $\Omega^*(\c{A})$ can be defined by 
relations of the form
\be
\theta^a \theta^b=P^{ab}{}_{cd} \theta^c \otimes \theta^d       \label{prod}
\ee
for suitable $P^{ab}{}_{cd} \in \c{Z}(\c{A})=\b{C}$. 
It can be shown that consistency with the nilpotency of $d$ requires
that the $\lambda_a$ satisfy a quadratic relation of the form
\be
2 \lambda_c \lambda_d P^{cd}{}_{ab} - 
\lambda_c F^c{}_{ab} - K_{ab} = 0.                   \label{manca}
\ee
The coefficients of the linear and constant terms must also belong to
the center. In the cases which interest us here they vanish. Notice
that Equation~(\ref{manca}) has the form of the structure equation of
a Lie algebra with a central extension.

We define~\cite{DubMadMasMou96b} the metric as a non-degenerate
$\c{A}$-bilinear map
\be
g: \Omega^1(\c{A}) \otimes_{\c{A}} \Omega^1(\c{A}) \rightarrow \c{A}.
\ee
This means that it can be completely determined up to central elements 
once its action on a basis of 1-forms is assigned. For example set
\be
g(\theta^a \otimes \theta^b)=g^{ab}.                          \label{metric}
\ee
The bilinearity implies that 
$$
fg^{ab}=g(f\theta^a \otimes \theta^b) =
g(\theta^a \otimes \theta^bf)=g^{ab}f
$$
and therefore $g^{ab} \in \c{Z}(\c{A})=\b{C}$: $\{\theta^a\}$ is a
special basis of 1-forms in which the coefficients of the metric are
central elements, namely complex numbers in our assumptions. This is
the property characterizing frames (vielbein) in ordinary geometry,
and is at the origin of the name `frame' for this basis also in
noncommutative geometry.  To define a covariant derivative $D$ which
satisfies~\cite{DubMadMasMou96b} a left and right Leibniz rule we
introduce a `generalized flip', an $\c{A}$-bilinear map
\be
\sigma: \Omega^1(\c{A}) \otimes_\c{A} \Omega^1(\c{A}) \rightarrow
\Omega^1(\c{A}) \otimes_\c{A} \Omega^1(\c{A}).             \label{2.2.4}
\ee
The flip is also completely determined once its action
on a basis of 1-forms is assigned. For example set
\be
\sigma (\theta^a \otimes \theta^b) = 
S^{ab}{}_{cd} \theta^c \otimes \theta^d.
\ee
As above, bilinearity implies that $S^{ab}{}_{cd}\in \c{Z}(\c{A})=\b{C}$.
Using the flip a left and right Leibniz rule can be written:
\bea
&& D (f \xi) =  df \otimes \xi + f D\xi  \\
&& D(\xi f) = \sigma (\xi \otimes df) + (D\xi) f.              
\eea
The torsion map 
\be
\Theta:\Omega^1(\c{A})\rightarrow\Omega^2(\c{A})
\ee
is defined by
\be
\Theta=d-\pi\circ D.
\ee
We shall assume that $\sigma$ satisfies the condition
\be
\pi \circ (\sigma + 1) = 0                                 \label{consis}
\ee
in order that the torsion be bilinear. The usual torsion 2-form
$\Theta^a$ is defined as $\Theta^a=d\theta^a-\pi\circ D\theta^a$.  It
is easy to check~\cite{DubMadMasMou96b} that if on the right-hand side
of Equation~(\ref{manca}) the term linear in $\lambda_a$ and the
constant term vanish then a torsion-free covariant derivative can be
defined by
\be
D \xi = - \theta \otimes \xi + 
\sigma (\xi \otimes \theta),                          \label{covdev0}
\ee
for any $\xi\in\Omega^1(\c{A})$.  The most general torsion-free $D$
for fixed $\sigma$ is of the form
\be
D = D_{(0)} + \chi                                           \label{2.2.15}
\ee
where $\chi$ is an arbitrary $\c{A}$-bimodule morphism
\be
\Omega^1(\c{A}) \buildrel \chi \over \longrightarrow
\Omega^1(\c{A}) \otimes \Omega^1(\c{A})                      \label{2.2.16}
\ee
fulfilling
\be
\pi \circ \chi = 0.                                            \label{2.2.20}
\ee
The compatibility of a covariant derivative with the metric is
expressed by the condition~\cite{DubMadMasMou95}
\be
g_{23}\circ D_2= d\circ g .                                \label{met-comp}
\ee
For the covariant derivative (\ref{covdev0}) this condition can be
written as the equation 
\be
S^{ae}{}_{df} g^{fg} S^{bc}{}_{eg} = g^{ab} \delta^c_d         \label{2.2.24}
\ee
if one uses the coefficients of the flip with respect to the frame.

Introduce the standard notation 
$\sigma_{12} = \sigma \otimes \mbox{id}$, 
$\sigma_{23} = \mbox{id} \otimes \sigma$, to extend to three factors
of a module any operator $\sigma$ defined on a tensor product of two
factors.  There is a natural continuation of the map (\ref{2.2.4}) to
the tensor product $\Omega^1(\c{A}) \otimes_\c{A} \Omega^1(\c{A})$
given by the map
\be
D_2(\xi \otimes \eta) = D\xi \otimes \eta + 
\sigma_{12} (\xi \otimes D\eta).                                \label{2.2.4e}
\ee
We define formally the curvature as
\be
\mbox{Curv} \equiv D^2= \pi_{12} \circ D_2 \circ D.                                              \label{curv}
\ee
We recover the standard definition of the frame components 
$R^a{}_{bcd}$ of the curvature tensor from the decomposition
\be
\mbox{Curv} (\theta^a) =
- {1 \over 2} R^a{}_{bcd} \theta^c \theta^d \otimes \theta^b    \label{2.16}
\ee
One can easily show \cite{FioMad99} that the curvature associated
to (\ref{covdev0}) is given by
\be
\mbox{Curv} (\xi )  = \xi_a \theta^2 \otimes \theta^a +
\pi_{12} \sigma_{12}\sigma_{23} \sigma_{12}
(\xi \otimes \theta \otimes \theta).                           \label{2.20}
\ee

The algebra we shall consider is a $*$-algebra. We shall require that 
the involution $*$ be extendable to the algebra of differential forms
in such a way that
\be
(\xi\eta)^* = (-1)^{pq}\eta^*\xi^*, \qquad \xi \in \Omega^p(\c{A}),
\quad \eta \in \Omega^q(\c{A}).                                  \label{sola}
\ee
We recall that the elements of the algebra are considered as 0-forms.
One would like to have a differential fulfilling the reality condition 
\be
(df)^* = df^*                                                 \label{d-real}
\ee 
as in the commutative case. Neither of the two differential
calculi we shall introduce in Section~\ref{preliminaries} satisfies
this condition; the differential calculus $\Omega^*(\c{A})$ is mapped
by $*$ into a new one $\bar\Omega^*(\c{A})$. As a consequence, the
reality conditions on the covariant derivative and curvature
formulated in \cite{FioMad98a} cannot be satisfied. However
we shall still suppose~\cite{DubMadMasMou95} that the extension of the
involution to the tensor product is given by
\be
(\xi \otimes \eta)^* = \sigma(\eta^* \otimes \xi^*).             \label{TPI}
\ee
A change in $\sigma$ therefore implies a change in the definition of
an hermitian tensor. The reality condition for the metric will be,
as in~\cite{FioMad98a},
\be
g\circ \sigma (\eta^* \otimes \xi^*) = 
(g(\xi \otimes \eta))^*.                                      \label{reality}
\ee 
We shall also continue to assume that
$\sigma$ satisfies the braid equation
\be
\sigma_{12}\sigma_{23} \sigma_{12}=\sigma_{23} \sigma_{12}\sigma_{23},
\ee
a condition implied~\cite{FioMad98a} by the reality condition on the
covariant derivative and the curvature. At the end of
Section~\ref{metrics} we shall briefly consider the question how to
modify reality condition on the covariant derivative and the curvature
in the present case.

\initiate
\section{The quantum Euclidean spaces and\\ 
their $q$-deformed differential calculi}               \label{preliminaries}

The starting point for the definition of the $N$-dimensional quantum
Euclidean space $\b{R}^N_q$ is the braid matrix $\hat{R}$ for
$SO_q(N,\b{C})$ a $N^2 \times N^2$ matrix, whose explicit expression
we give in Ap\-pen\-dix~\ref{append1}.  Certain properties of $\hat R$
which we shall use follow immediately from the definition. First, it
fulfills the braid equation
\be
\hat R_{12}\,\hat R_{23}\,\hat R_{12} =\hat R_{23}\,\hat R_{12}\,
\hat R_{23}.                                                  \label{braid1}
\ee
Here we have used again the conventional tensor notation
$\hat R_{12} = \hat R \otimes \mbox{id}$, 
$\hat R_{23} =  \mbox{id} \otimes \hat R $.
By repeated application of the Equation~(\ref{braid1}) one finds 
\be
f(\hat R_{12})\,\hat R_{23}\,\hat R_{12} 
=\hat R_{23}\,\hat R_{12}\,f(\hat R_{23})                     \label{braid2}
\ee
for any polynomial function $f(t)$ in one variable.  The
Equations~(\ref{braid1}) and~(\ref{braid2}) are evidently satisfied also
after the replacement $\hat R\rightarrow\hat R^{-1}$. Second, $\hat R$
is invariant under transposition of the indices:
\be
\hat R^{ij}_{kl}=\hat R^{kl}_{ij}.                             \label{propR1}
\ee
Here and in the sequel we use indices with values
\be
\begin{array}{ll}
i=-n,\ldots,-1,0,1,\ldots n, &
\mbox{ with $n\equiv \frac{N-1}{2}$ for $N$ odd}, \\[4pt]
i=-n,\ldots,-1,1,\ldots n,
& \mbox{ with $n \equiv \frac{N}{2}$ for $N$ even}.
\end{array}
\label{defn}
\ee
and $n$ to denote the rank of $SO(N,\b{C})$. The matrix element 
$\hat R^{ij}_{kl}$ vanishes unless the indices satisfy the following 
condition:
\be
\begin{array}{rl}
\mbox{either $i \neq -j$}  &\mbox{and $k=i$, $l=j$ or $l=i$, $k=j$}\\
\mbox{or $i=-j$}           &\mbox{and $k=-l$}.
\end{array}                                                  \label{propR3}
\ee
The $R$-matrix, defined by
$$
R^{ij}_{kl}=\hat R^{ji}_{kl},
$$
is lower-triangular.

There exists also~\cite{FadResTak89} a projector decomposition of 
$\hat R$:
\be
\hat R = q\c{P}_s - q^{-1}\c{P}_a + q^{1-N}\c{P}_t.       \label{projectorR}
\ee
The $\c{P}_s$, $\c{P}_a$, $\c{P}_t$ are $SO_q(N)$-covariant
$q$-deformations of the symmetric trace-free,
antisymmetric and trace projectors respectively and satisfy
the equations
\be
\c{P}_{\mu}\c{P}_{\nu} = \c{P}_\mu \delta_{\mu\nu}, \qquad
\sum_{\mu}\c{P}_{\mu} = 1, \quad \mu,\nu = s,a,t.
\ee
The $\c{P}_t$ projects on a one-dimensional sub-space and 
therefore it can be written in the form 
\be
\c{P}_t{}_{kl}^{ij} = (g^{sm}g_{sm})^{-1} g^{ij}g_{kl}
= \frac{k}{\omega_n (q^{-\rho_n+1}-q^{\rho_n-1})} g^{ij}g_{kl}  \label{Pt}
\ee
where $g_{ij}$ is the $N \times N$ matrix
\be
g_{ij}=q^{-\rho_i} \delta_{i,-j}.                               \label{defg}
\ee
We have here introduced the notation
$$
\rho_i=
\left\{
\begin{array}{ll}
(n-\frac{1}{2},\ldots,\frac{1}{2},0,-\frac{1}{2},\ldots,\frac{1}{2}-n)
& \mbox{ for $N$ odd}, \\[4pt]
(n-1,\ldots,0,0,\ldots,1-n) & \mbox{ for $N$ even.}
\end{array}
\right.
$$
and we have set
$$
k \equiv q-q^{-1}, \qquad
\omega_i \equiv q^{\rho_i}+q^{-\rho_i}.
$$
The matrix $g_{ij}$ is a $SO_q(N)$-isotropic tensor and is a 
deformation of the ordinary Euclidean metric
in a set of coordinates pairwise 
conjugated to each other under complex conjugation. 
It is easily verified  that its inverse $g^{ij}$ is given by
\be
g^{ij}=g_{ij}.
\ee
The metric and the braid matrix satisfy the 
relations~\cite{FadResTak89}
\be
g_{il}\,\hat R^{\pm 1}{}^{lh}_{jk} = 
\hat R^{\mp 1}{}^{hl}_{ij}\,g_{lk}, \qquad
g^{il}\,\hat R^{\pm 1}{}_{lh}^{jk} = 
\hat R^{\mp 1}{}_{hl}^{ij}\,g^{lk}.                          \label{gRrel}
\ee

The $N$-dimensional quantum Euclidean space is the associative algebra
$\b{R}^N_q$
generated by elements $\{x^i\}_{i=-n,\ldots,n}$ with
relations
\be
\c{P}_a{}^{ij}_{kl} x^k x^l=0.                              \label{xrel}
\ee 
These relations are preserved by the (right) 
action of the quantum group $U_qso(N)$,
which is defined on the generators by
\be
x^i\tl g=\rho_j^i(g) x^j,\qquad\qquad \rho_j^i(g)\in\b{C}, \label{trasf1}
\ee
where $\rho$ is the $N$-dimensional vector representation of $U_qso(N)$,
and extended to the rest of $\b{R}^N_q$  so that the latter becomes
a  $U_qso(N)$ module algebra. That is, for arbitrary $g,g'\in U_qso(N)$ 
and $a,a'\in\b{R}^N_q$, we have
\bea
a\tl (gg') &=& (a\tl g)\tl  g'                        \label{modalg1}\\
(aa')\tl g &=& (a\tl g_{(1)})\, (a'\tl g_{(2)}).       \label{modalg2}
\eea
Here we have used Sweedler notation (with lower indices) for the
coproduct, $\Delta(g)=g_{(1)}\otimes g_{(2)}$; the right-hand side 
is actually a short-hand notation for a finite
sum $\sum_Ig^I_{(1)}\otimes g^I_{(2)}$.

Relations (\ref{xrel}) can be written more explicitly in 
the form \cite{Ogi92a}
\be
\begin{array}{ll}
x^i x^j = q x^j x^i & \mbox{ for } i<j, i \neq -j, \\[3mm]
[x^i,x^{-i}]= k \omega_{i-1}^{-1} r^2_{i-1} 
&  \mbox{ for } i > 1, \\[4mm]
[x^1, x^{-1}]=\left\{
\begin{array}{l}
0 \\
h r_0^2
\end{array}
\right.&
\begin{array}{l}
\mbox{ for $N$ even,} \\
\mbox{ for $N$ odd.}
\end{array}
\end{array}
\label{explicitx}
\ee
We have here introduced $h$ defined by
\be
h \equiv q^{\frac{1}{2}}-q^{-\frac{1}{2}},
\ee
and we have defined as well a sequence of numbers $r_i$, $r$ by
\be
r^2_i=\sum_{k,l=-i}^i g_{kl} x^k x^l, \qquad   r^2\equiv r^2_n \label{defr}
\ee
where $i \ge 0$ for $N$ odd, whereas for $N$ even
$i\ge 1$ and of course in the sum only $k,l\neq 0$ actually occur.
The element $r^2$ is $SO_q(N)$-invariant and generates the center
of the algebra $\b{R}^N_q$.
It can be easily checked that
\be
x^j r_i^2=
\left\{
\begin{array}{lr}
r_i^2 x^j& \hbox{ for } |j| \le i, \nn[2pt]
q^2 r_i^2 x^j & \hbox{ for } j < -i, \\[2pt]
q^{-2} r_i^2 x^j & \hbox{ for } j > i. \nonumber
\end{array}
\right.
\label{xr2rel}
\ee

As this will be necessary for the construction of the elements $\lambda_a$
be introduced in section~4, we now extend the algebra 
$\b{R}^N_q$ by adding the square root $r_i$ of $r^2_i$ for $i=0 \ldots n$
as well as the inverses $r^{-1}_i$ of these elements. As the relations
(\ref{xr2rel}) contain only $q^{\pm 2}$ it is consistent to set
for $i \ge 0$
\be
x^j r_i=
\left\{
\begin{array}{lr}
r_i x^j& \mbox{ for } |j| \le i, \nn[2pt]
q r_i x^j & \mbox{ for } j < -i, \\[2pt]
q^{-1} r_i x^j & \mbox{ for } j > i. \nonumber
\end{array}
\right.                                                      \label{xrrel}
\ee
We shall be mainly interested in the case 
$q \in \b{R}^+$. In this case a conjugation
\be
(x^i)^*= x^j g_{ji}                                           \label{xconj}
\ee
can be defined on $\b{R}^N_q$ to obtain what is known as real quantum
Euclidean space. The elements $r_i$ are then real.

There are~\cite{CarSchWat91b} two differential calculi which are
covariant with respect to $U_qso(N)$, obtained by imposing the condition
\be
(d\alpha)\tl g=d (\alpha\tl g) \qquad\qquad\alpha\in\Omega^*(\b{R}_q^N)
\ee
on the differential. We denote the two exterior derivatives 
by $d$ and $\bar d$ and the corresponding exterior algebras by
$\Omega^*(\b{R}_q^N)$ and $\bar\Omega^*(\b{R}_q^N)$. 
If we introduce $\xi^i = dx^i$ and  $\bar{\xi}^i = \bar{d} x^i$, then 
they are characterized respectively by
\bea
x^i \xi^j &=& q\,\hat R^{ij}_{kl} \xi^k x^l,          \label{xxirel}\\[6pt]
x^i \bar \xi^j &=& 
q^{-1}\,\hat R^{-1}{}^{ij}_{kl} \bar\xi^k x^l.                \label{xxistar}
\eea
For $q\in\b{R}^+$ neither $\Omega^*(\b{R}_q^N)$ nor
$\bar\Omega^*(\b{R}_q^N)$ possesses an involution. However, one can
introduce a $*$-structure on the direct sum 
$\Omega^1(\b{R}_q^N) \oplus\bar\Omega^1(\b{R}_q^N)$ by setting
\be
(\xi^i)^* = \bar\xi^j g_{ji}.                               \label{xistar}
\ee
Using the properties (\ref{gRrel}, \ref{propR1}) of the 
$\hat R$-matrix one sees that the two calculi are conjugate; the
Equations (\ref{xxirel}) and (\ref{xxistar}) are exchanged.

By taking the differential of (\ref{xxirel}) and (\ref{xxistar}) the
$\xi \xi$-commutation relations are determined 
\be
\begin{array}{ll}
\c{P}_s{}_{kl}^{ij} \xi^k \xi^l = 0, &\qquad
\c{P}_t{}_{kl}^{ij} \xi^k \xi^l = 0,                   \label{xixirel} \\[6pt]
\c{P}_s{}_{kl}^{ij} \bar \xi^k \bar \xi^l = 0, &\qquad
\c{P}_t{}_{kl}^{ij} \bar \xi^k \bar \xi^l = 0. 
\end{array} \nonumber
\ee
These relations define the algebraic structure of $\Omega^*(\b{R}_q^N)$
and $\bar\Omega^*(\b{R}_q^N)$. 

It is useful to introduce a set of gradings $\deg_i$, $i=1,\ldots n$ on
$\Omega^*(\b{R}_q^N)$ by
\be
\deg_i(\xi^j)=\deg_i(x^j)=\left\{
\begin{array}{lr}
1 & \mbox{ if } i=j, \\ 
-1 & \mbox{ if } i=-j,\\
0 & \mbox{ otherwise.}
\end{array}
\right.                                                    \label{gradx}
\ee
All these gradings are preserved by the commutation
relations (\ref{xrel}), since the $\hat R$-matrix, and therefore any
polynomial function of it like $\c{P}_a$, fulfills (\ref{propR3}).
The $n$-ple $(\deg_1,...,\deg_n)$ coincides with the weight vector
of the fundamental vector representation of $so(N)$.

The Dirac operator \cite{Con94}, defined by Equation \ref{dirac},
\be
\xi^i = -[\theta,x^i]
\ee
is easily verified to be given by
\be
\theta= \omega_n q^{\frac{N}{2}} k^{-1} r^{-2} g_{ij} x^i \xi^j,
\ee
as pointed out in~\cite{Zum97, Ste96}.  For the barred calculus
$\bar\Omega^*(\b{R}_q^N)$ the `Dirac operator'
$\bar\theta$~(\ref{dirac}) is
\be
\bar\theta = -\omega_n q^{-\frac{N}{2}} k^{-1} 
r^{-2} g_{ij} x^i \bar\xi^j.                           \label{defbartheta}
\ee
If $q\in\b{R}^+$ it satisfies
\be
\theta^* = - \bar\theta.                                   \label{stern}
\ee

In order to construct the $\lambda_a$ and $\theta^a$ satisfying the
conditions described in Section~2 we first must solve the following
problem.  In Section~2 we assumed the center of the algebra $\c{A}$ to
be trivial, which makes possible the construction of elements
$\lambda_a$ and $\theta^a$ with the features described there.  But the
algebra generated by the $x^i$ and $r_j$ has a nontrivial center.
With a general Ansatz of the type
\be
\theta^a = \theta^a_i \xi^i                                \label{ansatz}
\ee 
the condition $[\theta^a,r_n^2]=0$ can be rewritten as
\be
(r^2_n \theta^a_i -q^{-2} \theta^a_i r^2_n) \xi^i=0,         \label{gaia}
\ee
which has no solution since $r^2_n\in\c{Z}(\b{R}_q^N)$.  To find a
solution to (\ref{gaia}) we further enlarge the algebra by adding a
unitary element $\Lambda$, the ``dilatator'', which satisfies the
commutation relations
\be
x^i \Lambda=q \Lambda x^i.                                 \label{xLambda}
\ee

We also add its inverse $\Lambda^{-1}$.  In the case $N$ odd we can
now follow the scheme previously proposed for $N=3$~\cite{FioMad99}
but in the case of even $N$ the situation is slightly more
complicated. We have added the elements 
$r_1^{\pm 1}=(x^{-1} x^1)^{\pm \frac{1}{2}}$ and as a consequence the
center is non trivial even after addition of $\Lambda$.  The elements
$$
r_1^{-1}x^{\pm 1} = \left(x^1(x^{-1})^{-1}\right)^{\pm \frac{1}{2}}
$$
commute also with $\Lambda$. (We recall that $x^{-1})^{-1}$ is the
inverse of the element $x^i$ with $i = -1$.)
In other words, since the algebra generated by  $(\Lambda, r_1^{\pm 1},
x^{\pm 1}, \ldots x^{\pm n})$ is completely symmetric in the exchange of
$x^1$ and $x^{-1}$, there is no way to distinguish between these two elements.
To have $N$ linearly independent $\theta^a$, instead of fewer, we shall 
need to add yet another element to the algebra.
We choose to add a `Drinfeld-Jimbo' 
generator $K=q^{\frac {H_1}2}$ and its inverse $K^{-1}$, 
where $H_1$ belongs to the Cartan subalgebra of $U_qso(N)$ and represents 
the component of the angular momentum in the $(-1,1)$-plane. 
This new element satisfies the commutation relations
\be
\begin{array}{lll}
K x^{\pm 1}= q^{\pm 1} x^{\pm 1} K,   &&\\
K x^{\pm i}= x^{\pm i} K,             &&\mbox{for $i>1$},
\end{array}                                               \label{xkapparel}
\ee
as well as
\be
K \Lambda=\Lambda K.                                   \label{kappaLambda}
\ee
When $q\in\b{R}^+$ it is compatible with the commutation relations to
extend the $*$-structure (\ref{xconj}) $\Lambda,K$ as
\be
\Lambda^*=\Lambda^{-1}, \qquad\qquad K^*=K.
\ee
We must decide now which commutation relations $\Lambda,K$ should
satisfy with the $\xi^i$. As already observed \cite{FioMad99} there 
are different possibilities. 

A first possibility is to set~\cite{Ogi92a} 
\be
\xi^i \Lambda=\Lambda \xi^i, \qquad \Lambda d=q d \Lambda. \label{xiLambda}
\ee
This choice has the disadvantage that $\Lambda$ cannot be considered
as an element of the quantum space, because due to (\ref{xiLambda}) it
does not satisfy the Leibniz rule 
$d(fg)=f dg + (df) g \;\forall f,g \in \b{R}^N_q$.  Nevertheless, it can
be interpreted in a consistent way as an element of the Heisenberg
algebra, because $\Lambda^{-2}$ can be constructed~\cite{Ogi92a} as a
simple polynomial in the coordinates and derivatives.

Alternatively, what was considered
also in \cite{AscCastSca97}, one could ask the Leibniz rule 
$d(fg)=f dg + (df) g $ to hold also if $f=\Lambda$.
By differentiating (\ref{xLambda}) one obtains that
\be
\xi^i \Lambda+x^i d \Lambda = 
q d \Lambda x^i + q \Lambda \xi^i.                       \label{lambdanec}
\ee
A solution would be to require that
\be
x^i (d \Lambda)=q (d \Lambda) x^i, \qquad 
\xi^i \Lambda=q \Lambda \xi^i.                         \label{altxiLambda}
\ee
In particular it would then be possible to set $d\Lambda = \Lambda d$,
which implies that $(d \Lambda)=0$.  This choice is not completely
satisfactory either since we would like the relation
\be
df=0 \mbox{ implies } f \propto 1                            \label{center}
\ee
to hold, and this would not be the case if $d \Lambda=0$. As a
consequence the general formalism is still not strictly applicable and
there will be a conformal ambiguity in the choice of metric. We shall
see below that with a procedure similar to the one described
previously \cite{FioMad99} for $N=3$, we would recover 
$\b{R} \times S^{N-1}$ as geometry rather than $\b{R}^N$ in the
commutative limit.  Therefore, in the sequel we shall impose the first
condition (\ref{xiLambda}). As will be shown in the next section, this
allows us to normalize the $\theta^a$ and $\lambda_a$ in such a way as
to obtain $\b{R}^N$ as geometry in the commutative limit.
The above discussion with $\Lambda$ can be repeated to determine the
commutation relations between $K$ and the 1-forms $\xi^i$.  We choose
$d K=0$. Then consistency with (\ref{xkapparel}) requires that
\be
\begin{array}{lll}
K \xi^{\pm 1}=q^{\pm 1} \xi^{\pm 1} K, &&\\
K \xi^i=\xi^i K,                       &&\mbox{for $i>1$},
\end{array}                                               \label{xikapparel}
\ee

To summarize, we shall consider the algebra $\c{A}_N$, an extension of 
$\b{R}^N_q$ defined for odd $N$ as 
$$
\c{A}_N = \{x^i,\;r_j,\;r_j^{-1},\;\Lambda,\Lambda^{-1} 
\st -n \le i \le n,\; 0 \le j<n\}
$$
with generators which satisfy the relations (\ref{xrel}), 
(\ref{xrrel}), (\ref{xLambda}) and for even $N$ as 
$$
\c{A}_N = \{x^i,\;r_j,\;r_j^{-1},\;\Lambda,\; \Lambda^{-1},\;  
K,\;K^{-1} \st -n \le i \le n,\; 1 \le j<n\}
$$
with generators which satisfy the relations (\ref{xrel}),
(\ref{xrrel}), (\ref{xLambda}), (\ref{xkapparel}).  The algebra of
differential forms $\Omega^*(\c{A}_N)$ is generated by the one-forms
$\xi^i$ satisfying relations (\ref{xxirel}), (\ref{xixirel}),
(\ref{xiLambda}) when $N$ is odd, and (\ref{xxirel}), (\ref{xixirel}),
(\ref{xiLambda}), (\ref{xikapparel}) when $N$ is even.  However one
must bear in mind that the additional elements $\Lambda$ and $K$ are
rather exceptional since $dK=0$ and either $d \Lambda=0$, or it does
not satisfy the Leibniz rule. These elements would be better
interpreted as elements of the Heisenberg algebra.

\initiate
\section{Inner derivations and frame}

We would like to construct a frame $\theta^a$ and the associated inner
derivations $e_a = \ad \lambda_a$ satisfying the conditions in
Section~2 for the case of the algebra $\c{A}_N$. We first solve this
problem in a larger algebra, which we now define.  It is possible to
extend $\Omega^*(\c{A}_N)$ to the cross-product algebra
\mbox{$\Omega^*(\c{A}_N)\cocross\uqs$} by postulating the 
cross-commutation relations
\be
\xi g = g_{(1)} (\xi\tl\:g_{(2)} )                           \label{gxrel}
\ee
for any $g\in\uqs$ and $\xi\in\Omega^*(\c{A}_N)$. The algebra
$\Omega^*(\c{A}_N)\cocross\uqs$ can be made into a module algebra
under the action $\tl$ of $\uqs$ by extending the latter on the
elements of $\uqs$ as the adjoint action,
$$
h\tl g=Sg_{(1)} h g_{(2)},          \qquad g,h\in \uqs.
$$
The $S$ here denotes the antipode of $\uqs$.

Let us introduce $\uqso$ the Hopf algebra with the same algebra 
structure of $\uqs$,  but opposite coalgebra, and by 
$\Delta^{op}(g)=g_{(2)}\otimes g_{(1)}$
its coproduct. On any module algebra $\c{M}$ of $\uqso$ the 
corresponding action $\tlo$ will thus fulfill the relations
\bea
a\tlo(gg') &=& (a\tlo g) \tlo   g'                      \label{modalg1o}\\
(aa')\tlo g &=& (a\tlo g_{(2)})\, (a'\tlo g_{(1)}).       \label{modalg2o}
\eea
These are to be compared with (\ref{modalg1}) and (\ref{modalg2}).
It is immediate to show that definition (\ref{gxrel}) implies  
that one can realize the action $\tl$ in the `adjoint-like way'
\be
\eta\tl g=Sg_{(1)} \, \eta \, g_{(2)}                 \label{realiz1}
\ee
on all of $\Omega^*(\c{A}_N)\cocross\uqs$. 
On the other hand, one can realize also a corresponding action $\tlo$ by
\be
\eta\tlo g=( S^{-1}g_{(2)}) \, \eta \, g_{(1)},                      
                                                       \label{realiz2}
\ee
where $S^{-1}$ is the antipode of $\Delta^{op}$.

We return now to the problem of the construction of a frame and
of a set of dual inner derivations for the differential calculus
$(\Omega^*(\c{A}_N), d)$. As a first step, we must find $N$ independent
solutions $\vartheta^a$ to the equation
\be
[f,\vartheta^a]=0 \qquad\forall f\in\c{A}_N.            \label{framecond}
\ee
We shall look first for solutions $\vartheta^a$ in 
$\Omega^*(\c{A}_N)\cocross\uqs$. The reason is the following.
For each solution $\vartheta^a$ of (\ref{framecond}) 
and for any $g\in\uqs$ we can consider the image
$\vartheta_g{}^a\in \Omega^*(\c{A}_N)\cocross\uqs$ of $\vartheta^a$, defined by
\be
\vartheta_g{}^a:=\vartheta^a\tlo g=(S^{-1}g_{(2)})\,\vartheta^a\,g_{(1)}.
                                                       \label{frametransf}
\ee
Now we show that its commutator with any element $f\in\c{A}_N$
vanishes:
\bea
[f,\vartheta_g{}^a] &=& [f, (S^{-1}g)_{(1)}\,\vartheta^a\,S(S^{-1}g)_{(2)}]\nn
 &=&f\,(S^{-1}g)_{(1)}\,\vartheta^a\,S(S^{-1}g)_{(2)}-
(S^{-1}g)_{(1)}\,\vartheta^a\,S(S^{-1}g)_{(2)}\,f\nn
&\stackrel{(\ref{gxrel})}{=}
&(S^{-1}g)_{(1)}\,[f\tl (S^{-1}g)_{(2)}]\vartheta^a\,S(S^{-1}g)_{(3)} \nn
&& - (S^{-1}g)_{(1)}\,\vartheta^a\,[f\tl (S^{-1}g)_{(2)}]\,S(S^{-1}g)_{(3)}\nn
&=&
(S^{-1}g)_{(1)}\,\left[f \tl (S^{-1}g)_{(2)},
\vartheta^a\right]S(S^{-1}g)_{(3)}=0
\nonumber
\eea
In the last equality we have used the fact that $f\tl (S^{-1}g)_{(2)}
\in\c{A}_N$ and (\ref{framecond}).  In other words $\vartheta_g{}^a$
also behaves as a frame element. Moreover by its very definition
$\vartheta_g{}^a$ will in general belong to
$\Omega^1(\c{A}_N)\cocross\uqs$ even if
$\vartheta^a\in\Omega^1(\c{A}_N)$ (unless the $n$-plet
$\{\vartheta^a\}$ builds an irreducible representation of $\uqso$).
We can summarize the above results in the
\begin{prop}
\label{framesubspace}
The subspace ${\cal F}$
of $\Omega^1(\c{A}_N)\cocross\uqs$ spanned by frame elements carries
a representation of $\uqso$.
\end{prop}
We shall call a set $\{\vartheta^a\}_{a=-n,...,n}$ a generalized frame
if $\vartheta^a$ are elements of $\Omega^1(\c{A}_N)\cocross\uqs$
fulfilling (\ref{framecond}), and any $\xi \in \Omega^1(\c{A}_N)\cocross\uqs$ 
can be uniquely decomposed in the form $\vartheta^a\xi_a$ with 
$\xi_a\in\c{A}_N\cocross\uqs$.

We now look for a generalized frame in the form of a basis of an
irreducible $N$-dimensional representation of $\uqso$. Recall that the
universal $R$-matrix $\R$ is a special element
\be
\R\equiv\R^{(1)}\otimes\R^{(2)}\in\uqs\otimes \uqs
\label{swee}
\ee
intertwining between  $\Delta$ and
$\Delta^{op}$, and so does also $\R^{-1}_{21}$:
\be
\R\Delta(\cdot)=\Delta^{op}(\cdot)\R, \qquad\qquad
\R^{-1}_{21}\Delta(\cdot)=\Delta^{op}(\cdot)\R^{-1}_{21}.  \label{inter}
\ee
In  (\ref{swee}) we have used a Sweedler notation
with upper
indices: the right-hand side is a short-hand notation
for a sum $\sum_I\R_I^{(1)}\otimes\R_I^{(2)}$ of infinitely 
many terms. The other main properties of $\R$ are recalled in 
the appendix \ref{UnivR}. 
\begin{prop}
Let $(u_5)^{-1}=\R^{(1)}(S\R^{(2)})$,
$(u_7)^{-1}=\R^{-1}{}^{(2)}(S\R^{-1}{}^{(1)})$.  
The elements of $\Omega^1(\c{A}_N)\cocross\uqs$
\bea
&&\vartheta^a:=\alpha \:u_7^{-1}\left(S\R^{(2)}\right) 
\xi^a\,\R^{(1)}  \label{def2}\\
&&\bar \vartheta^a :=\bar\alpha \:u_5^{-1}\left(S\R^{-1}{}^{(1)}
\right)\xi^a\, \R^{-1}{}^{(2)}                  \label{def2bar}
\eea
are covariant under the action (\ref{realiz2}), more precisely
\be
\vartheta^a\tlo g=\rho_b^a(g)\vartheta^b \qquad\qquad
\bar\vartheta^a\tlo g=\rho_b^a(g)\bar \vartheta^b .         \label{opact2}
\ee
\end{prop}
In the previous definitions we have inserted two scalar factors
$\alpha,\bar\alpha$ to be fixed later.

\bp{}
We prove the first formula. We find
\bea
\vartheta^a\tlo g & \stackrel{(\ref{realiz2})}{=}& (S^{-1}g_{(2)})
\vartheta^a
g_{(1)}\stackrel{(\ref{def2})}{=}(S^{-1}g_{(2)})\alpha \:u_7^{-1}
\left(S\R^{(2)}\right) \xi^a\,\R^{(1)} g_{(1)}           \nonumber\\
& \stackrel{(\ref{interantip})}{=}& \alpha\: u_7^{-1}(Sg_{(2)})
\left(S\R^{(2)}\right) \xi^a\,\R^{(1)} g_{(1)}           \nonumber\\
& \stackrel{(\ref{inter})}{=}& \alpha\: u_7^{-1}
\left(S\R^{(2)}\right)(Sg_{(1)}) \xi^a\, g_{(2)}\R^{(1)} \nonumber\\
& \stackrel{(\ref{trasf1})}{=}& 
\rho_b^a(g)\alpha\: u_7^{-1}\left(S\R^{(2)}\right) \xi^b\, \R^{(1)}
\stackrel{(\ref{def2})}{=}\rho_b^a(g)\vartheta^b.         \nonumber
\eea
The proof of the second formula is completely analogous.
\ep
We can  give an alternative and very useful expression for 
$\vartheta^a,\bar\vartheta^a$. It is convenient to introduce the
Faddeev-Reshetikin-Taktadjan generators \cite{FadResTak89} of $\uqs$:
\be
\c{L}^+{}_l^a:=\R^{(1)}\rho_l^a(\R^{(2)})\qquad\qquad
\c{L}^-{}_l^a:=\rho_l^a(\R^{-1}{}^{(1)})\R^{-1}{}^{(2)}. \label{frt}
\ee
In our conventions  
\be
\R\in\uqsp\otimes\uqsn
\ee
where $\uqsp,\uqsn$ denote the positive and negative Borel subalgebras;
hence we see that 
$\c{L}^+{}_l^a\in\uqsp$ and $\c{L}^-{}_l^a\in\uqsn$. 
{From} formulae (\ref{delta1}), (\ref{delta2}) in the 
Appendix~\ref{UnivR} one finds that the coproducts are given by
\be
\Delta(\c{L}^+{}^i_j)=\c{L}^+{}^i_h\otimes\c{L}^+{}^h_j \qquad\qquad
\Delta(\c{L}^-{}^i_j)=\c{L}^-{}^i_h\otimes\c{L}^-{}^h_j.
\label{coprodL}
\ee
Their commutation relations will be given in following section. 
Using them and (\ref{gxrel}) one easily proves the
\begin{prop} 
Let $u_3=\R^{(2)}S^{-1}\R^{(1)}$, 
$u_4=\R^{-1}{}^{(1)}S^{-1}\R^{-1}{}^{(2)}$.
\bea
&&\vartheta^a=\alpha\c{L}^-{}^a_l\eta^l= \alpha 
\eta^m  \rho_m^l(u_4) \c{L}^-{}^a_l          \label{utile1}\\
&&\bar \vartheta^a=\bar\alpha\c{L}^+{}^a_l\eta^l=\bar\alpha 
\eta^m\rho_m^l(u_3)\c{L}^+{}^a_l,        \label{utile2}
\eea
with $\eta^i=\xi^i$ or $\eta^i=\bar\xi^i$.
Thus $\vartheta^a$ and $\bar \vartheta^a$ belong to 
$\Omega^1(\c{A}_N) \cocross \uqsn$ and
$\Omega^1(\c{A}_N) \cocross \uqsp$ respectively for $\eta^i=\xi^i$, or
$\bar\Omega^1(\c{A}_N) \cocross \uqsn$ and
$\bar\Omega^1(\c{A}_N) \cocross \uqsp$ respectively for $\eta^i=\bar\xi^i$.
\end{prop}
We now show that $(\vartheta^a,\bar \vartheta^a)$ constitute a frame.
Recall that the braid matrix can be written as
$\hat R^{ij}_{hk}=(\rho^j_h\otimes\rho^i_k)\R$. Using
(\ref{gxrel}), (\ref{coprodL}) one can easily prove
\begin{lemma} 
\bea
x^i\c{L}^{\pm}{}^a_b&=&\c{L}^{\pm}{}^a_c x^j\hat 
R^{\pm 1}{}^{ci}_{jb},\\
\xi^i\c{L}^{\pm}{}^a_b&=&\c{L}^{\pm}{}^a_c 
\xi^j\hat R^{\pm 1}{}^{ci}_{jb},                          \label{lemmino}\\
\bar\xi^i\c{L}^{\pm}{}^a_b&=&\c{L}^{\pm}{}^a_c 
\bar\xi^j\hat R^{\pm 1}{}^{ci}_{jb}
\eea
\end{lemma}
We should note that the commutation relations we have just found
are different from (in a sense opposite to)
the  ones which we would have found by imposing, instead of (\ref{gxrel}),
the condition
\be
g \xi=g_{(1)}\tr \xi g_{(2)},                            \label{alter}
\ee
with a {\it left} action $\tr$. This is true also in the
commutative case ($q=1$), and in a sense is unpleasent
since the latter definition
of the commutation relations is what we are usually more familiar to. 
For instance, if $h$ is a primitive generator
in the Cartan subalgebra then $[h,x^i]$ defined in the first way
is opposite to the one defined in the second; or if $g^+$ is a positive
root, then $[g^+,x^i]$ defined in the first way is proportional (up to
a Cartan subalgebra factor)
to the commutator $[g^-,x^i]$ defined in the second way, where 
$g^-$ is the negative root opposite to $g^+$.

The reason why we have adopted (\ref{gxrel}) is that this is necessary
in order that the coordinates $x^i$ carry {\it upper} 
indices (as it is conventional
in general relativity) and that the
representation $\rho$ defined by (\ref{trasf1}) can be considered as the
fundamental (vector) one, rather than its contragradient. This follows
from $\rho^i_h(gg')=\rho^i_j(g)\rho^j_h(g')$. Had we used lower indices
to label the coordinates,
or replaced $\rho(g)$ with $\rho(Sg)$, we could have adopted
(\ref{alter}).

\begin{prop}
Assume $\xi^i,\bar\xi^i$ are the differentials (\ref{xxirel}),
(\ref{xxistar}), of the previous section. If we choose
\be
\bar\alpha=\alpha^{-1}=\Lambda                    \label{choice}
\ee 
then $\{\vartheta^a\}$, $\{\bar\vartheta^a\}$ are 
generalized frames respectively in 
$\Omega^1(\c{A}_N)\cocross\uqs$
and $\bar\Omega^1(\c{A}_N) \cocross \uqs$.
\end{prop}
\bp{}
We prove this in the first case.
$$
x^i\vartheta^a\stackrel{(\ref{utile1})}{=}x^i\c{L}^-{}^a_l\Lambda^{-1}\xi^l 
\stackrel{(\ref{lemmino})}{=}\c{L}^-{}^a_c \hat R^{-1}{}^{ci}_{jl}
x^j\Lambda^{-1}\xi^l \stackrel{(\ref{xxirel})}{=}
\c{L}^-{}^a_c\Lambda^{-1}\xi^c x^i
\stackrel{(\ref{lemmino})}{=}\vartheta^a x^i
$$
The proof in the second case is completely analogous.
\ep

In the commutative limit $\{\xi^i\}_{i=1,...,n}$ is a frame and all
other frames can be obtained from it by a $G\subset Uso(N)$
transformation.  Formulae (\ref{utile1}) and (\ref{utile2}) say that in
the noncommutative case one frame can be obtained from $\{\xi^i\}$ by
a very particular $\uqs$ transformation, the one with matrix elements
$\c{L}^{\pm}{}^a_b$.

Note that by the choice (\ref{choice}) $\vartheta^a$, $\bar\vartheta^a$
commute not only with the coordinates $x^i,r_j$, but also with the special
element $\Lambda$,
\be
[\Lambda, \vartheta^a]=0 \qquad\qquad [\Lambda, \bar\vartheta^a]=0,
\ee
under the assumption (\ref{xiLambda}); had we assumed instead
(\ref{altxiLambda}), then the same would be true by choosing
$\alpha=\Lambda^{-1}r$, $\bar\alpha=\Lambda r$.
On the other hand, if we choose 
$\bar\alpha=\bar f(r)\Lambda$, $\alpha=f(r)\Lambda^{-1}$
(with $f,\bar f\neq$const), then $\vartheta^a$, $\bar\vartheta^a$ will still
commute with the coordinates $x^i,r_j$, but in general not with $\Lambda$.

The commutation relations between the frame elements are given by the
\begin{prop}
The commutation relations among the $\vartheta^a$ (resp.
$\bar\vartheta^a$) are  as the ones among the $\xi^i$ (resp. 
$\bar\xi^i$), except for the opposite products:
\be
\begin{array}{ll}
\c{P}_s{}_{ab}^{cd}\vartheta^b\vartheta^a=0,  
\qquad &\c{P}_t{}_{ab}^{cd}\vartheta^b\vartheta^a=0       \cr
 \c{P}_s{}_{ab}^{cd}\bar\vartheta^b\bar\vartheta^a=0,  
\qquad &\c{P}_t{}_{ab}^{cd}\bar\vartheta^b\bar\vartheta^a=0.   
\end{array}                               \label{ththrel}
\ee
\end{prop}
\bp{}
We prove the claim in the first case (in the second the proof
is completely analogous). For both $\c{P}=\c{P}_s,\c{P}_t$
\bea
&&\c{P}_{ab}^{cd}\vartheta^b\vartheta^a
\stackrel{(\ref{utile1})}{=}\c{P}^{ab}_{cd}\Lambda^{-1}\c{L}^-{}^b_l\xi^l
\Lambda^{-1}\c{L}^-{}^a_m\xi^m\stackrel{(\ref{lemmino})}{=}
\Lambda^{-2}\c{P}_{ab}^{cd}\c{L}^-{}^b_l\c{L}^-{}^a_h \xi^k
\hat R^{-1}{}^{hl}_{mk}\xi^m
\nonumber\\
&&\quad\stackrel{(\ref{L+L+rel})}{=}\Lambda^{-2}\xi^l\xi^n 
(\c{P}\hat R^{-1})_{ln}^{mh}\c{L}^-{}_h^d  \c{L}^-{}_m^c
\stackrel{(\ref{projectorR})}{\propto}\Lambda^{-2}\xi^l\xi^n 
\c{P}_{ln}^{mh} \c{L}^-{}_h^d  \c{L}^-{}_m^c\stackrel{(\ref{xixirel})}{=}0.
\nonumber
\eea
\ep
We now decompose $d$ (and similarly $\bar d$) in terms both of
differentials $\xi^i=dx^i$ and $\uqs$-covariant derivatives
$\partial_i$ and of frame elements $\vartheta^a$ and derivations
$\epsilon_a$:
\be
d=\partial_i\xi^i=\epsilon_a\vartheta^a \qquad\qquad
d=\bar\partial_i\bar\xi^i=\bar \epsilon_a\bar\vartheta^a.
\ee
Looking at (\ref{utile1}), (\ref{utile2}) we find that
\be
\partial_i=\epsilon_a \c{L}^-{}^a_i\Lambda^{-1}\qquad\qquad
\bar\partial_i=\bar \epsilon_a \c{L}^+{}^a_i\Lambda.
                                              \label{concile}
\ee

Now assume that the ``Dirac operator'' exists; it must necessarily be 
$\uqs$-invariant. We can decompose it as well  both in the basis of 
differentials and in the basis of frame elements:
\be
\theta=-w_i\xi^i=-y_a\vartheta^a,  \qquad\qquad
\bar\theta=-\bar w_i\bar \xi^i=-\bar y_a\bar \vartheta^a,     \label{dec} 
\ee
with $w_i,\bar w_i\in\c{A}_N$ and $y_a\in\c{A}_N\cocross\uqsn$,
$\bar y_a\in\c{A}_N\cocross\uqsp$; in this case
\be
\epsilon_a=[y_a,\cdot]\qquad\qquad\bar \epsilon_a=[\bar y_a,\cdot].
\ee
The commutation relations among the $w_i, \bar w_i$ will be of the form
\be
\c{P}_a{}_{ij}^{hk}w_kw_h=0, \qquad
\c{P}_a{}_{ij}^{hk}\bar w_k\bar w_h=0.            \label{wwrel}
\ee
{From} the \uqs invariance of $d$, $\bar d$, $\theta$ and $\bar\theta$ 
it follows that
\be
\mu_i\tl g=\mu_l\,\rho_i^l(Sg),   \qquad\qquad \mu_i=w_i,\bar w_i,
\partial_i,\bar\partial_i.
\ee
{From} (\ref{dec}), (\ref{utile1}), (\ref{utile2})
we find
\bea
&&y_a=w_i S\c{L}^-{}^i_a  \Lambda                \label{yaexp}\\
&&\bar y_a=\bar w_i S\c{L}^+{}^i_a  \Lambda^{-1}  \label{byaexp}
\eea
{From} the $\uqso$-invariance of $d,\bar d,\theta,\bar\theta$ we find the
transformation rules
\be
\nu_a\tlo g=\nu_b\,\rho^b_a(S^{-1}g) , \qquad\qquad\nu_a=y_a,\bar y_a, 
\epsilon_a,\bar \epsilon_a.                                      \label{OPact}
\ee
\begin{prop}
The commutation relations among the $y_a$ and among the $\bar y_a$ 
are 
\be
\c{P}_a{}_{ab}^{cd}y_cy_d=0, \qquad\qquad
\c{P}_a{}_{ab}^{cd}\bar y_c\bar y_d=0.             \label{yyrel}
\ee
\end{prop}
\bp{}
\bea
\c{P}_a{}_{ab}^{cd}y_cy_d &\propto & \c{P}_a{}_{ab}^{cd}\,
w_i (S\c{L}^-{}^i_c) w_j (S\c{L}^-{}^j_d) \Lambda^2\nn
&=&\c{P}_a{}_{ab}^{cd}\,w_i (w_j\tl \c{L}^-{}^i_h)(S\c{L}^-{}^h_c)  
(S\c{L}^-{}^j_d)  \Lambda^2\nn 
&=&\c{P}_a{}_{ab}^{cd}\,w_i w_k\rho^k_j(S\c{L}^-{}^i_h)(S\c{L}^-{}^h_c)  
(S\c{L}^-{}^j_d)  \Lambda^2\nn
&=& w_i w_k \hat R^{ki}_{hj}S\left(\c{L}^-{}^j_d
\c{L}^-{}^h_c\c{P}_a{}_{ab}^{cd}\right)  \Lambda^2\nn
&\stackrel{(\ref{L+L+rel}),(\ref{frt})}{=}
&w_i w_k \hat R^{ki}_{hj}S\left(\c{P}_a{}_{lm}^{hj}\c{L}^-{}^m_b
\c{L}^-{}^l_a\right)  \Lambda^2 \nn
&\propto& w_i w_k \c{P}_a{}_{lm}^{ki} S\left(\c{L}^-{}^m_b
\c{L}^-{}^l_a\right)  \Lambda^2 
\stackrel{(\ref{wwrel})}{=}0. \nonumber
\eea
The proof is completely analogous for $\bar y_a$.
\ep
{From} $\xi^i=[x^i,\theta]=[w_j\xi^j,x^i]$,  
$\bar\xi^i=[x^i,\bar\theta]=[\bar w_j\bar\xi^j,x^i]$ it follows
$$
x^iw_k\Lambda=\hat R^{-1}{}^{ji}_{hk}w_j\Lambda\!x^h-
\delta^i_h\Lambda\qquad\qquad
x^i\bar w_k\Lambda^{-1}=\hat R^{ji}_{hk}\bar w_j\bar 
\Lambda\!x^h-\delta^i_h\Lambda^{-1}
$$
using these relations and
(\ref{gxrel}), (\ref{coprodL}) one can easily prove
\begin{prop}
\be
[y_a,x^j]= \Lambda S\c{L}^-{}^j_a
\qquad\qquad
[\bar y_a,x^j]= \Lambda^{-1} S\c{L}^+{}^j_a \label{yxrel}
\ee
\end{prop}
These formulae can be seen as the inverse of (\ref{yaexp}), (\ref{byaexp}),
since they give $S\c{L}^+{}^a_j,S\c{L}^-{}^a_j$ in terms of $y_a,\bar y_a$. 

In next section we find algebra homomorphisms
\bea
&&\varphi^+:\c{A}_N\cocross \uqsp \rightarrow\c{A}_N,\\
&&\varphi^-: \c{A}_N\cocross\uqsn \rightarrow\c{A}_N,
\eea
defined as the identity on $\c{A}_N$, 
\be
\varphi^{\pm}(a)=a \qquad\mbox{ if }a\in\c{A}_N.
\ee
Then the 1-forms
\bea
&&\theta^a=\theta_l^a\xi^l\qquad\qquad 
\theta_l^a:=\varphi^-(\Lambda^{-1}\c{L}^-{}^a_l)       \label{deftheta}\\
&&\bar\theta^a=\bar\theta_l^a\bar\xi^l\qquad\qquad 
\bar\theta_l^a:=\varphi^+(\Lambda \c{L}^+{}^a_l)
\eea
belong to $\Omega^1(\c{A}_N),\bar\Omega^1(\c{A}_N)$ and still fulfil
the property (\ref{framecond}), in other words build up a frame.  
Similarly the elements of $\c{A}_N$ defined by
\be
\lambda_a=\varphi^-(y_a) \qquad\qquad\bar\lambda_a=\varphi^+(\bar y_a).
\ee
will yield the dual inner derivations,
\be
e_a = [\lambda_a,\cdot]\qquad\qquad 
\bar e_a = [\bar \lambda_a,\cdot]
\ee
It is clear that since $\c{Z}(\c{A}_N)=\b{C}$ then the frame and the
dual set of inner derivations are {\it uniquely} determined up to a
linear transformation (with numerical coefficients).  By definition,
the $\lambda_a,\bar\lambda_a$ will still fulfill the commutation
relations (\ref{yyrel}),
\be
\c{P}_{(a)}{}^{ab}_{cd} \lambda_a \lambda_b=0, \qquad\qquad
\c{P}_{(a)}{}^{ab}_{cd} \bar\lambda_a \bar \lambda_b=0. \label{lalarel}
\ee

After application of $\varphi^{\pm}$ equations (\ref{concile}) become
\be
\partial_i=e_a \varphi^-(\c{L}^-{}^a_i)\Lambda^{-1}\qquad\qquad
\bar\partial_i=\bar e_a \varphi^+(\c{L}^+{}^a_i)\Lambda.
                                              \label{pconcile}
\ee
Their interest lies in the fact that they relate the
$SO_q(N)$-covariant derivatives $\partial_i,\bar\partial_i$
\cite{CarSchWat91b,OgiSchWesZum92}, introduced following the approach
of Woronowicz and Wess-Zumino~\cite{Wor79,Wor87a,Wor87b,WesZum90}, to
the inner derivations $e_a$, $\bar e_a$ dual to the frame and which
are defined using ordinary commutation relations, following the
approach of Connes~\cite{Con86,Con94}.

On the other hand, $\varphi^+,\varphi^-$ cannot be extended to
homomorphisms of 
$\Omega^*(\c{A}_N)\cocross U^{\pm}_qso(N)\rightarrow\Omega^*(\c{A}_N)$, 
since the commutation relations between the elements of $\uqs$ and the
1-forms do not map into the commutations of the elements of
$\c{A}_N$ with the 1-forms; this is immediate to check if we use the
frame basis: the frame elements don't commute with $U^{\pm}_qso(N)$,
but do commute with
$\varphi^+(\uqsp),\varphi^-(\uqsn)\subset\c{A}_N$. As a consequence,
the commutation relations among the $\theta_a,\bar\theta_a$ differ
from (\ref{ththrel}) by a reversal of the product order. That is,
\begin{prop}
The commutation relations among the $\theta^a$ (resp. $\bar\theta^a$)
are  as the ones among the $\xi^i$ (resp. $\bar\xi^i$):
\be
\begin{array}{ll}
\c{P}_s{}_{ab}^{cd}\theta^a\theta^b=0, \qquad\qquad  
\c{P}_t{}_{ab}^{cd}\theta^a\theta^b=0 \cr            
\c{P}_s{}_{ab}^{cd}\bar\theta^a\bar\theta^b=0, \qquad\qquad  
\c{P}_t{}_{ab}^{cd}\bar\theta^a\bar\theta^b=0 \cr            
\end{array}                                               \label{thetherel}
\ee
\end{prop}
\bp{}
We prove the claim in the first case. For both $\c{P}=\c{P}_s,\c{P}_t$
we have
\bea
\c{P}_{ab}^{cd}\theta^a\theta^b &\stackrel{(\ref{deftheta})}{=}&
\c{P}_{ab}^{cd}\theta^a\theta^b_n\xi^n
\stackrel{(\ref{framecond})}{=}   
\c{P}_{ab}^{cd}\theta^b_n\theta^a\xi^n             \nonumber\\
&\stackrel{(\ref{deftheta})}{=}&\c{P}_{ab}^{cd}\Lambda^{-2}
\varphi^-(\c{L}^-{}_n^b)\varphi^-(\c{L}^-{}_l^a)\xi^l\xi^n
\propto\Lambda^{-2} \varphi^-(\c{P}_{ab}^{cd}\c{L}^-{}_n^b
\c{L}^-{}_l^a)\xi^l\xi^n                         \nonumber\\
&\stackrel{(\ref{L+L+rel})}{=}&\Lambda^{-2}\varphi^-(\c{L}^-{}_b^d
\c{L}^-{}_a^c\c{P}^{ab}_{ln})
\xi^l\xi^n\stackrel{(\ref{xixirel})}{=}0\nonumber
\eea
The proof is completely analogous for $\bar y_a$.
\ep

The $\theta^a,\bar\theta^a,\lambda_a,\bar\lambda_a$ do not inherit from
$\vartheta^a,\bar\vartheta^a,y_a,\bar y_a$ the transformation properties
(\ref{opact2}), (\ref{OPact}) under the action $\tlo$ of \uqso.
For the $\theta^a,\bar\theta^a$  this is clear since the
commutation relations (\ref{thetherel}) are no longer
compatible with the action of \uqso. As a consequence, this is true also
for the $\lambda_a,\bar\lambda_a$, because 
$\theta=-\theta^a\lambda_a$, $\bar\theta=-\bar\theta^a\bar\lambda_a$ 
are invariant. 

One could in principle define an \uqso action $\tlo'$ on $\c{A}_N$ 
by postulating instead
$$
\lambda_a\tlo'g =\lambda_b\rho^b_a(S^{-1}g)  \qquad\qquad
\bar \lambda_a\tlo' g=\bar \lambda_b\rho^b_a(S^{-1}g) ;
$$
but the latter would differ from the action $\tlo$ 
fulfilling (\ref{realiz2}). We shall therefore not do so.

\initiate
\section{Homomorphism $\c{A}_N\cocross U_qso(N)\rightarrow \c{A}_N$} 

It is known \cite{FadResTak89} that a set of generators of 
$U_qso(N)$ is provided by \linebreak
$\{\c{L}^+{}^i_j,\c{L}^-{}^i_j\}$ and some further elements
obtained by introducing square roots and inverses of the
elements $\c{L}^{\pm}{}^i_i$, which is always possible
because they are invertible elements belonging to the Cartan subalgebra. 
The set $\{\c{L}^+{}^i_j,\c{L}^-{}^i_j\}$ has $2N^2$ 
elements, but $N(N-1$) of them vanish due to the
upper and lower triangularity of the matrices $\c{L}^{\pm}$, 
whereas their diagonal elements are the inverses of each other:
\bea
&&\c{L}^+{}^i_j=0,           \hspace{1.5cm}\mbox{if $i>j$}\label{sfilza1}\\
&&\c{L}^-{}^i_j=0,           \hspace{1.5cm}\mbox{if $i<j$}\label{sfilza2}\\
&&\c{L}^{+}{}^i_i\c{L}^{-}{}^i_i=1,\hspace{1cm}\forall i  \label{sfilza3}\\
&&\c{L}^{\pm}{}^{-n}_{-n}....\c{L}^{\pm}{}^n_n=1.         \label{sfilza4}
\eea
To these relations we have to add the relations characterizing 
$U_qso(N)$, 
\be
\c{L}^{\pm}{}^i_j\c{L}^{\pm}{}^h_k g^{kj}=g^{hi}\qquad 
\c{L}^{\pm}{}_i^j\c{L}^{\pm}{}_h^k g_{kj}=g_{hi},  \label{LLg}
\ee
which further reduce to $N(N-1)/2$ the number of independent 
generators, and finally the commutation relations
\bea
&&\hat R^{ab}_{cd}\,\c{L}^+{}^d_f\c{L}^+{}^c_e=
\c{L}^+{}^b_c\c{L}^+{}^a_d\,\hat R^{dc}_{ef}      \label{L+L+rel}\\
&&\hat R^{ab}_{cd}\,\c{L}^-{}^d_f\c{L}^-{}^c_e=
\c{L}^-{}^b_c\c{L}^-{}^a_d\,\hat R^{dc}_{ef} \label{L-L-rel}\\
&&\hat R^{ab}_{cd}\,\c{L}^+{}^d_f\c{L}^-{}^c_e=
\c{L}^-{}^b_c\c{L}^+{}^a_d\,\hat R^{dc}_{ef}.\label{L+L-rel}
\eea
Note that (\ref{sfilza4}) yields no constraint for $N$ even since it is a 
consequence of (\ref{LLg}), whereas for $N$ odd it yields one
further constraint. In fact (\ref{LLg}) implies that 
$(\c{L}^+{}^0_0)^2=1$ which together with (\ref{sfilza4})
implies that $\c{L}^+{}^0_0=1$.
The antipode on the generators takes the form 
\be
S\c{L}^{\pm}{}^i_j=g^{hi}\c{L}^{\pm}{}^k_hg_{jk}       \label{antip}
\ee

To construct a homomorphism 
$\varphi:\c{A}_N\cocross U_qso(N)\rightarrow \c{A}_N$ acting as the
identity on $\c{A}_N$ it is therefore sufficient to define it on
$\c{L}^+{}^i_j,\c{L}^-{}^i_j$ and to verify that all of relations
(\ref{lemmino}) and (\ref{sfilza1}) to (\ref{L+L-rel}) are satisfied.
Applying $\varphi$ to (\ref{yxrel}) and using (\ref{antip}) it follows
that
\bea
&&\varphi(\c{L}^-{}^i_j)\stackrel{(\ref{antip})}{=} 
g^{ih}\varphi(S\c{L}^-{}^k_h)g_{kj}\stackrel{(\ref{yxrel})}{=}
g^{ih}\Lambda^{-1}[\lambda_h,x^k]g_{kj} \label{lambdaxrel}\\
&&\varphi(\c{L}^+{}^i_j)\stackrel{(\ref{antip})}{=} 
g^{ih}\varphi(S\c{L}^+{}^k_h)g_{kj}\stackrel{(\ref{yxrel})}{=}
g^{ih}\Lambda[\bar\lambda_h,x^k]g_{kj}.          \label{barlambdaxrel}
\eea
The problem is thus equivalent to the construction of $N$ objects
$\lambda_a$ and $N$ objects $\bar\lambda_a$ such that relations
(\ref{lambdaxrel}), (\ref{barlambdaxrel})
define elements $\varphi(\c{L}^{\pm}{}^k_h)$
in $\c{A}_N$ which satisfy all of relations (\ref{lemmino})$_1$  and
(\ref{sfilza1}) to (\ref{L+L-rel}).

Actually for the construction of a frame and a set of dual inner
derivations in $\Omega^*(\c{A}_N)$ (resp.  $\bar\Omega^*(\c{A}_N)$) we
just need a homomorphism $\varphi^-:U_q^-so(N)\cross\c{A}_N\rightarrow
\c{A}_N$ (resp. $\varphi^+:U_q^+so(N)\cross\c{A}_N\rightarrow
\c{A}_N$) acting as the identity on $\c{A}_N$, what we look for first.
This is equivalent to looking for $N$ objects $\lambda_a$ (resp.
$\bar\lambda_a$) such that the objects
$\varphi^-(\c{L}^-{}^i_j)=g^{ih}\Lambda^{-1}[\lambda_h,x^k]g_{kj}$
(resp. $\varphi^+(\c{L}^+{}^i_j)=g^{ih}\Lambda[\bar\lambda_h,x^k]g_{kj}$)
fulfill just the relations involving only $\c{L}^-$ (resp. $\c{L}^+$).
Finally we find which further conditions $\lambda_a,\bar\lambda_a$
must fulfill in order that we can `glue' $\varphi^-,\varphi^+$ into a
unique homomorphism $\varphi$.

\begin{theorem}                                         \label{theor1}
One can define a homomorphism 
$\varphi^-:\c{A}_N\cocross U^-_qso(N)\rightarrow \c{A}_N$ by
setting on the generators
\bea
&& \varphi^-(a)=a,         \qquad \forall a\in\c{A}_N, \\[4pt]
&& \varphi^-(\c{L}^-{}^i_j)=g^{ih}\Lambda^{-1}[\lambda_h,x^k]g_{kj},
\eea
with
\be
\begin{array}{ll}
\lambda_0=\gamma_0 \Lambda (x^0)^{-1}
&\quad\mbox{for $N$ odd,} \\[6pt]
\lambda_{\pm 1}=\gamma_{\pm 1} \Lambda (x^{\pm 1})^{-1} K^{\mp 1} 
&\quad\mbox{for $N$ even,} \\[6pt]
\lambda_a=\gamma_a \Lambda r_{|a|}^{-1}r_{|a|-1}^{-1} x^{-a}
&\quad\mbox{otherwise,} 
\end{array}                                             \label{deflambda}
\ee
and $\gamma_a \in \b{C}$ normalization constants fulfilling
the conditions
\be
\begin{array}{ll}
\gamma_0 = -q^{-\frac{1}{2}} h^{-1} &\quad\mbox{for $N$ odd,} \\[6pt]
\gamma_1 \gamma_{-1}=
\left\{\begin{array}{l}
-q^{-1} h^{-2}\\
k^{-2}
\end{array}\right.
&\quad\!\begin{array}{l}
\mbox{for $N$ odd,} \\
\mbox{for $N$ even,}
\end{array}\\[8pt]
\gamma_a \gamma_{-a} =
-q^{-1} k^{-2} \omega_a \omega_{a-1} &\quad\mbox{for $a>1$}. \nonumber
\end{array}                                               \label{gamma}
\ee
\end{theorem}
The proof is given in Appendix \ref{proofs}.  The relation
(\ref{gamma}) fixes only the product $\gamma_a \gamma_{-a}$.  Notice
that $\gamma_0^2$ for $N$ odd and $\gamma_1 \gamma_{-1}$ for $N$ even
are positive real numbers, while all the remaining products 
$\gamma_a \gamma_{-a}$ are negative.  Notice that the embedding
$\c{A}_N \hookrightarrow \c{A}_{N+2}$ defined by (\ref{embedding})
automatically induces an embedding for the corresponding $\lambda_a$.
Similarly one can prove
\begin{theorem}                                         \label{theor2}
One can define a homomorphism
$\c{A}_N\cocross U^+_qso(N)\rightarrow \c{A}_N$ by
setting on the generators
\bea
&& \varphi^+(a)=a,           \qquad \forall a\in\c{A}_N, \\[4pt]
&& \varphi^+(\c{L}^+{}^i_j)=g^{ih}\Lambda[\bar\lambda_h,x^k]g_{kj},
\eea
with
\be
\begin{array}{ll}
\bar\lambda_0=\bar\gamma_0 \Lambda^{-1} (x^0)^{-1}
&\quad\mbox{for $N$ odd,} \\[6pt]
\bar\lambda_{\pm 1} = 
\bar\gamma_{\pm 1} \Lambda^{-1} (x^{\pm 1})^{-1} K^{\pm 1} 
&\quad\mbox{for $N$ even,} \\[6pt]
\bar\lambda_a = 
\bar\gamma_a \Lambda^{-1} r_{|a|}^{-1}r_{|a|-1}^{-1} x^{-a}
&\quad\mbox{otherwise,} 
\end{array}                                         \label{defbarlambda}
\ee
and $\bar\gamma_a \in \b{C}$ normalization constants fulfilling
the conditions
\be
\begin{array}{ll}
\bar \gamma_0 = q^{\frac{1}{2}} h^{-1}
&\quad\mbox{for $N$ odd,} \\[6pt]
\bar \gamma_1 \bar \gamma_{-1} =
\left\{
\begin{array}{l}
-q h^{-2}\\
k^{-2}
\end{array}
\right.
&\quad\!\begin{array}{l}
\mbox{for $N$ odd,} \\
\mbox{for $N$ even,}
\end{array}\\[8pt]
\bar \gamma_a \bar \gamma_{-a}=
-q k^{-2} \omega_a \omega_{a-1}
&\quad\mbox{for $a>1$}.         \nonumber
\end{array}                                             \label{bargamma}
\ee
\end{theorem}
To `glue' $\varphi^+,\varphi^-$ into a unique homomorphism 
$\varphi: \c{A}_N\cocross U_qso(N)\rightarrow \c{A}_N$ we still need
to satisfy the image under $\varphi$ of equations (\ref{L+L-rel}) and
(\ref{sfilza3}).  As we shall prove in appendix \ref{proofs2}, this is
possible only in the case of odd $N$ and completely fixes the
coefficients $\gamma_a,\bar\gamma_a$ in the previous Ans\"atze:
\begin{theorem}                                          \label{theor3}
In the case of odd $N$ we one can define a homomorphism
$\uqs \cross\c{A}_N\rightarrow \c{A}_N$
by setting on the generators
\bea
&& \varphi(a)=a\qquad \forall a\in\c{A}_N \\[4pt]
&& \varphi(\c{L}^-{}^i_j)=g^{ih}\Lambda^{-1}[\lambda_h,x^k]g_{kj},\\[4pt]
&& \varphi(\c{L}^+{}^i_j)=g^{ih}\Lambda[\bar\lambda_h,x^k]g_{kj},
\eea
with $\lambda_j,\bar\lambda_j$ defined as in (\ref{deflambda}),
(\ref{defbarlambda}) and with coefficients given by
\be
\begin{array}{lr}
\gamma_0 = -q^{-\frac 12}h^{-1}\\[4pt]
\gamma_1^2 = -q^{-2} h^{-2} \\[4pt]
\gamma_a^2 = -q^{-2}\omega_a \omega_{a-1} k^{-2} 
&\quad \mbox{for $a>1$}\\[4pt]
\gamma_a = q \gamma_{-a}  
& \quad \mbox{for $a \le 1$} \\[4pt]
\bar \gamma_a=-q \gamma_a.
\end{array}                                              \label{fixgamma} 
\ee
\end{theorem}
Notice that the $\gamma_a,\bar\gamma_a$ for $a \neq 0$ are imaginary
and fixed only up to a sign. This has as a consequence that the
homomorphism $\varphi$ does not preserve the star structure of \uqs,
in other words it is not a star-homomorphism.
In the case $N=3$ this is the same homomorphism which is also
constructed in \cite{CerMadSchWes00}.
Alternatively, we can fix the normalizations of the barred objects
so that for $q\in\b{R}^+$ the involution gives
\be
\lambda_a^* = -g^{ab} \bar\lambda_b, \qquad 
(\theta^a)^*= \bar\theta^b g_{ba}.                    \label{invlambda}
\ee
Then the coefficients $\bar\gamma_a$ 
will be related to the $\gamma_a$ by
\be
\begin{array}{lll}
\bar\gamma_0=-q\gamma_0^*, &&\mbox{if $N$ odd,}\\[4pt]
\bar\gamma_{\pm 1}=-\gamma_{\mp 1}^*, 
&&\mbox{if $N$ even,}                    \label{def*lambda}\\[4pt]
\bar\gamma_a=-\gamma_{-a}^*
\left\{
\begin{array}{ll}
1   &\;\mbox{if $a>0$}\\
q^2 &\;\mbox{if $a<0$}
\end{array}
\right.
&&\mbox{otherwise.}
\end{array}
\ee

We summarize our results for the frames $\theta^a,\bar\theta^a$:
\bea
&&\theta^a=\theta^a_l\xi^l=
\Lambda^{-2}g^{ab}[\lambda_b,x^j]g_{jl}\xi^l \\[6pt]
&&\bar\theta^a=\bar\theta^a_l\bar\xi^l=
\Lambda^2g^{ab}[\bar\lambda_b,x^j]g_{jl}\bar\xi^l.
\eea
These objects commute both with $x^i, r_j$ and with $\Lambda$.
Had we adopted the commutation rule (\ref{altxiLambda}),
instead of (\ref{xiLambda}), then the same would be true by 
introducing an additional factor $r$ in the right-hand side 
of both the previous equations. The matrix elements 
$\bar\theta^a_i,\theta^a_i$ fulfill the $\varphi^+,\varphi^-$ 
images of equation (\ref{L+L+rel}),  (\ref{L-L-rel})
\be
\hat R_{cd}^{ab} \theta^d_j \theta^c_i = \theta^b_l
\theta^a_k\hat R_{ij}^{kl}\qquad\qquad
\hat R_{cd}^{ab} \bar\theta^d_j \bar\theta^c_i = \bar\theta^b_l
\bar\theta^a_k\hat R_{ij}^{kl}.
\label{rtheta} 
\ee

\initiate
\section{Metrics and linear connections on the \\ 
quantum Euclidean space}                                \label{metrics}

In this Section we construct the covariant derivative, the 
corresponding linear connection and the metric associated with the 
frames introduced in Section~5. We follow for $\b{R}^N_q$ the same 
scheme proposed previously~\cite{FioMad99} for $\b{R}^3_q$.

Since the $\Omega^1(\c{A}_N)$ and $\bar \Omega^1(\c{A}_N)$, we recall, 
are free modules, the covariant derivatives
can be defined by their actions on the frame
\bea
D\theta^a &=& - \omega^a{}_{bc} \theta^b \otimes \theta^c, \\
\bar D \bar \theta^a &=& - \bar \omega^a{}_{bc} \bar \theta^b
\otimes \bar \theta^c.         \nonumber
\eea
For the generalized permutation $\sigma$ we can write 
\be
\sigma (\theta^a\otimes \theta^b) =
S^{ab}_{cd}\,\theta^c\otimes \theta^d.
\label{sigma}
\ee
For the same reasons as for the coefficients of the metric the
requirement of bilinearity
\be
f S^{ab}_{cd} \theta^c \otimes \theta^d = 
\sigma (f \theta^a \otimes \theta^b) = 
\sigma (\theta^a \otimes \theta^b f) =
S^{ab}_{cd} f \theta^c \otimes \theta^d
\ee
forces $S^{ab}{}_{cd} \in \c{Z}(\c{A}_N)=\b{C}$.
According to the consistency condition (\ref{consis})
for the torsion and the commutation relations (\ref{thetherel})
for the frame we find that the matrix $S$ has the form 
\be
S = C_s\c{P}_s - \c{P}_a + C_t\c{P}_t,
\label{defS}
\ee
with $C_s$ and $C_t$ complex $N^2\times N^2$ matrices.
As the next step we would like to define the metric according to
(\ref{metric}).
\be
g(\theta^a \otimes \theta^b)=g^{ab}.
\label{metric2}
\ee
But we have a problem. Due to the form of (\ref{projectorR})
it is not possible to satisfy simultaneously the metric compatibility
condition and (\ref{defS}). Similarly to what was previously
done for $N=3$ \cite{FioMad99} the best we can do is
to weaken the compatibility condition to a
condition of proportionality.
In this way like in the case $N=3$ we find the two solutions
\be
S = q\hat R, \qquad\qquad S = (q\hat R)^{-1},       \label{double}
\ee
corresponding to the choices $C_s = q^2$, $C_t=q^{2-N}$ or 
$C_s = q^{-2}$, $C_t =q^{N-2}$ respectively. As $\hat R$ does, also both
these solutions for $S$  satify the Yang-Baxter equation (\ref{braid1}).
Now, using the property (\ref{gRrel}) for the $\hat R$-matrix, we see that
\be
S^{ae}_{df} g^{fg} S^{cb}_{eg}=q^{\pm 2} g^{ac} \delta^b_d.
\ee
This has to be compared with the metric compatibility condition
(\ref{2.2.24}), 
which in a basis becomes
\be
S^{ae}_{df} g^{fg} S^{cb}_{eg}=g^{ac} \delta^b_d.
\ee
The metric is in fact compatible with the linear connection only up to a
conformal factor.

We can compute the action of $\sigma$ on the basis $\xi^i \otimes \xi^j$
\be
\sigma(\xi^i \otimes \xi^j)=S^{ij}_{hk} \xi^h \xi^k.
\label{sigmaxi}
\ee
To see this we have started with the definition (\ref{sigma}) of the action
on $\theta^a \otimes \theta^b$, and used (\ref{ansatz}).
Then the result follows from (\ref{rtheta}).
Notice that (\ref{sigmaxi}) coincides with (\ref{sigma}).
In a similar way, from (\ref{metric2}), (\ref{ansatz}) and (\ref{gtt})
we can determine the action of $g$ on the $\xi$.
\be
g(\xi^i \otimes \xi^j)= g^{ij} \Lambda^2.
\label{sigmag}
\nonumber
\ee
This expression contains $\Lambda$, therefore the addition of this element
to the algebra is a necessary condition for the construction of the metric,
even if we would perform the calculation directly in the basis given
by $\xi^i$.
According to (\ref{covdev0}) a covariant derivative can be defined by
\be
D \xi = 
-\theta \otimes \xi + \sigma (\xi \otimes \theta).
\label{covder}
\ee
As $\theta$ is $SO_q(N)$-invariant, so is $D$. This is true for both
choices of $\sigma$. One can show~\cite{FioMad99} that the 
expression~(\ref{covder}) is the most general torsion-free, 
$SO_q(N)$-invariant linear connection.

The explicit action of the covariant derivative on $\xi^i$
can be computed. We apply (\ref{gRrel}), (\ref{xxirel}) and obtain for
$S=q^{-1} \hat R^{-1}$ 
\be
D\xi^i=0.
\ee
Analogously, from (\ref{gRrel}), (\ref{xxirel}), (\ref{squareR}) it turns 
out that
\bea
D\xi^i &=&-\theta \otimes \xi^i+q^{2+\frac{N}{2}} \omega_n k^{-1}
r^{-2} g_{hj} (\hat R^2)^{ji}_{lm} x^h \xi^l \otimes \xi^m \nn[4pt] 
&=& (q^2-1) \theta \otimes \xi^i - q^2 \omega_n r^{-2} g_{lm}
\left( q^{1-\frac{N}{2}} x^i \xi^l \otimes \xi^m
+q^{\frac{N}{2}} \hat R^{mi}_{hj} x^l \xi^h\otimes \xi^j \right)
\nonumber\\[4pt]
&=& (q^2-1) (\theta \otimes \xi^i + \xi^i\otimes \theta)
-q^{3-\frac{N}{2}} \omega_n r^{-2} x^i\xi^l\otimes \xi^mg_{lm}.
\label{explicitD} 
\eea
The coordinates are well adapted to the first linear connection, but not to
the second one, because in the latter case $D\xi^i\neq 0$.

However, for both of these covariant derivatives the corresponding 
curvature (\ref{curv}) vanishes.
\be
\mbox{Curv}(\xi)=0
\ee
This can be seen by performing the same calculation done previously
\cite{FioMad99} for $N=3$.
\begin{eqnarray*}
\mbox{Curv}(\xi) && = \pi_{12}\sigma_0{}_{12}
\sigma_0{}_{23}\sigma_0{}_{12}
(\theta^a\otimes \theta^b\otimes \theta^c)\:\xi_a\lambda_b\lambda_c \\[4pt]
&& = \pi_{12}(S_{12}S_{23}S_{12})^{abc}{}_{def}
(\theta^d\otimes \theta^e\otimes \theta^f)\:\xi_a\lambda_b\lambda_c \\[4pt]
&& = (S_{12}S_{23}S_{12})^{abc}{}_{def}
(\theta^d\theta^e\otimes \theta^f)\:\xi_a\lambda_b\lambda_c \\[4pt]
&& = -(S_{12}S_{23}\c{P}_a{}_{12})^{abc}{}_{def}
(\theta^d\theta^e\otimes \theta^f)\:\xi_a\lambda_b\lambda_c \\[4pt]
&& =  -(\c{P}_a{}_{23}S_{12}S_{23})^{abc}{}_{def}
(\theta^d\theta^e\otimes \theta^f)\:\xi_a\lambda_b\lambda_c \\[4pt]
&& = 0. 
\end{eqnarray*}
We have used here the braid relation (\ref{braid1}) for $S$, (\ref{consis})
and (\ref{lalarel}).

To conclude this section we repeat the above construction for the
calculus $\bar\Omega^1(\c{A}_N)$.  Again, there is a unique
metric $g$ and two torsion-free $SO_q(N)$-covariant linear connections
compatible with it up to a conformal factor.  In the $\bar \theta^a$
basis the actions of $g$ and $\sigma$ are respectively
\bea
&&g(\bar\theta^a\otimes\bar\theta^b) = g^{ab} \\[2pt]
&&\sigma (\bar\theta^a\otimes \bar\theta^b) =
\bar S^{ab}{}_{cd}\,\bar\theta^c\otimes \bar\theta^d,
\eea
The two choices for $\sigma$ are
\be
\bar S = q\hat R, \qquad\mbox{or}\qquad \bar S = (q\hat R)^{-1}.
\ee
This implies
\be
\bar S {}^{ae}{}_{df} g^{fg} \bar S {}^{cb}{}_{eg} 
= q^{\pm 2} g^{ac} \delta^b_d.
\ee
In the $\bar\xi^i$ basis the actions of $g$ and $\sigma$ become
\bea
&&g(\bar \xi^i\otimes \bar\xi^j) = 
g^{ij}  \Lambda^{-2}, \\[2pt]
&&\sigma (\bar \xi^i\otimes \bar \xi^j) =
\bar S^{ij}{}_{hk}\,\bar\xi^h\otimes \bar\xi^k.
\eea
The two covariant derivatives, one for each choice of $\sigma$, are
\be
\bar D \bar \xi = -\bar \theta \otimes \bar \xi + 
\sigma (\bar \xi \otimes \bar \theta).
\ee
The associated linear curvatures $\overline{\mbox{Curv}}$ vanish.

Since in the commutative limit $q \rightarrow 1$ we have
$$
\Lambda \rightarrow 1, \qquad K \rightarrow 1,
$$
we have also
$$
g(\xi^i\otimes \xi^j)\rightarrow\delta^{i,-j} \qquad
g(\bar \xi^i\otimes \bar \xi^j)\rightarrow\delta^{i,-j}.
$$
The right-hand side is the matrix of coefficients of the flat metric
in the complex cartesian coordinates $x^i$, and we recover $\b{R}^N$
as geometry (at least formally). Had we adopted the commutation rule
(\ref{altxiLambda}), instead of (\ref{xiLambda}), then we would have
found an additional factor $\lim_{q\to 1}r^2$ in the right-hand side,
which would have corresponded to the coefficients of the metric of
$\b{R}\times S^{N-1}$.

\initiate
\section{Appendix}
\subsection{Miscellanea}                                 \label{append1}

In this appendix we give some miscellanea formulae
on $\b{R}_q^N$.

The braid matrix of $SO_q(N)$ is given by \cite{FadResTak89}
\bea
\hat R&=&q \sum_{i \neq 0} \delta^i_i \otimes \delta^i_i +
\sum_{\stackrel{\scriptstyle i \neq j,-j} 
{\mbox{ or } i=j=0}} \delta^j_i \otimes \delta^i_j+ q^{-1} 
\sum_{i \neq 0} \delta^{-i}_i
\otimes \delta^i_{-i} 
\label{defR} \\
&&+k (\sum_{i<j} \delta^i_i \otimes \delta^j_j- 
\sum_{i<j} q^{-\rho_i+\rho_j} 
\delta^{-j}_i \otimes \delta^j_{-i}) \nonumber
\eea
where here $\delta^i_j$ is the $N \times N$ matrix with all elements
equal to zero except for a $1$ in the $i$th column and $j$th row.
Clearly, 
$R^{ij}_{hk}:=\hat R^{ji}_{hk}$
is a lower-triangular matrix.
Moreover, under substitution of $q$ by $q^{-1}$ in (\ref{defR}) we find that
\be
\hat R(q^{-1}){}^{ij}_{kl}=\hat R(q)^{-1} {}_{-k-l}^{-i-j}.    \label{propR2}
\ee
Using the projector decomposition (\ref{projectorR}) and the expression
(\ref{Pt}) for $\c{P}_t$ the square of the $\hat R$-matrix can be computed
to be
\be
(\hat R^2)^{ij}_{kl}=k \hat R^{ij}_{kl}+\delta^i_k \delta^j_l 
-q^{1-N} k g^{ij} g_{kl}.                                  \label{squareR}
\ee

There is an embedding $\c{A}_N \hookrightarrow \c{A}_{N+2}$
given by
\bea
&x^i \rightarrow x^i &\quad\mbox{for}\quad -n \le i \le n,  \nn
&r_i \rightarrow r_i &\quad\mbox{for}\quad 0 \le i \le n
\label{embedding} \\[2pt]        
&\Lambda \rightarrow \Lambda,& K \rightarrow K.\nonumber
\eea

It follows from (\ref{explicitx}) that one can rewrite (\ref{defr}) as
\be
r^2_i =\omega_i \left(q^{-1} \omega_{i-1}^{-1} r^2_{i-1}+ x^i x^{-i}
\right) =\omega_i \left(q \omega_{i-1}^{-1}r^2_{i-1}+ x^{-i} x^i\right),
\label{rutil1}
\ee
for $i>1$ and for $N$ odd and $i=1$. This implies that for $i \ge 1$ 
\be
x^{-i} x^i-q^2 x^i x^{-i}=-q k \omega_i^{-1} r^2_i.
\label{rutil2}
\ee
Finally, we recall a useful property of the fundamental representation
$\rho$ of $U_qso(N)$, namely
\be
\rho^a_b(Sg)=g^{ad}\rho^c_d(g)g_{cb}.
\ee

\subsection{Universal $R$-matrix}                          \label{UnivR} 

In this appendix we recall the basics about the universal $R$-matrix
\cite{Dri86} of a quantum group $\uqg$, while fixing our conventions.
We recall some useful formulae 
\bea
&&(\Delta \otimes \mbox{id})\R=\R_{13}\R_{23} \label{delta1}\\
&&(\mbox{id}\otimes\Delta  )\R=\R_{13}\R_{12} \label{delta2}\\
&&(S \otimes \mbox{id})\R=\R^{-1}= (\mbox{id}\otimes S^{-1})\R\\
&&S^{-1}(g) = u^{-1}S(g) u.                             \label{interantip}
\eea
Here $u$ is any of the elements $u_1,u_2,..u_8$ defined below:
\be
\begin{array}{ll}
u_1:= (S\R^{(2)}) \R^{(1)} \qquad\qquad &u_2:= (S\R^{-1}{}^{(1)})
\R^{-1}{}^{(2)} \cr
u_3:= \R^{(2)} S^{-1}\R^{(1)} \qquad\qquad &u_4:= \R^{-1}{}^{(1)} 
S^{-1}\R^{-1}{}^{(2)} \cr
(u_5)^{-1}:= \R^{(1)} S\R^{(2)} \qquad\qquad &(u_6)^{-1}:= 
(S^{-1}\R^{(1)}) \R^{(2)} \cr
(u_7)^{-1}:= \R^{-1}{}^{(2)} S\R^{-1}{}^{(1)} 
\qquad\qquad &(u_8)^{-1}:= 
(S^{-1}\R^{-1}{}^{(2)}) \R^{-1}{}^{(1)} 
\label{defu}
\end{array}
\ee
In fact, using the results of Drinfel'd \cite{Dri86,Dri90} one can show that
\be
u_1=u_3=u_7=u_8=vu_2=vu_4=vu_5=vu_6,
\ee
where $v$ is a suitable element belonging to the center of \uqs.

\medskip

{From} (\ref{inter}) and (\ref{delta1},\ref{delta2}) it follows the
universal Yang-Baxter relation
\be
\R_{12}\R_{13}\R_{23}=\R_{23}\R_{13}\R_{12},          \label{YBEQ1}
\ee
whence the other two relations follow
\bea
\R^{-1}{}_{12}\R^{-1}{}_{13}\R^{-1}{}_{23}& =& 
\R^{-1}{}_{23}\R^{-1}{}_{13}\R^{-1}{}_{12}\label{YBEQ2}\\
\R_{13}\R_{23}\R^{-1}{}_{12}& =& 
\R^{-1}{}_{12}\R_{23}\R_{13}\label{YBEQ3}
\eea
By applying $\mbox{id}\otimes \rho^a_c\otimes\rho^b_d$ to 
(\ref{YBEQ1}), $\rho^a_c\otimes\rho^b_d\otimes\mbox{id}$ to
(\ref{YBEQ2}) and $\rho^a_c\otimes\mbox{id}\otimes\rho^b_d$
to (\ref{YBEQ3}) we respectively find the commutation relations
(\ref{L+L+rel}), (\ref{L-L-rel}), (\ref{L+L-rel}).

\subsection{Proof of Theorem \ref{theor1}}                    \label{proofs}

We divide the proof in various steps. First note that, because of
(\ref{gRrel}), equation (\ref{lemmino}) for
$\varphi^-(\c{L}^-)$ is equivalent to
\be
x^h [\lambda_a,x^i]=q\hat{R}^{hi}_{jk}[\lambda_a, x^j] x^k. 
\label{Rxlambda}
\ee
\begin{prop}
\label{prop1}
$N$ independent solutions to Equation~(\ref{Rxlambda}) are given by
(\ref{deflambda})
where $\gamma_a \in \b{C}$ are arbitrary normalization constants.
\end{prop}
\bp{}
As a first step in the proof that
Equation~(\ref{Rxlambda}) is satisfied by the $\lambda_a$ of
(\ref{deflambda}), we calculate the commutation relations between the
$x^i$ and the $\lambda_a$:
\be
\begin{array}{lll}
\Lambda\varphi(\c{L}^-{}^{-a}_{-i})
=q^{\rho_a-\rho_i}[\lambda_a, x^i]&=0& \mbox{ for $i<a$,} \\[4pt]
\Lambda\varphi(\c{L}^-{}^{-a}_{-i})
=q^{\rho_a-\rho_i}[\lambda_a, x^i]&=-q^{\rho_a-\rho_i+1} k \lambda_a x^i
& \mbox{ for $i>a$,} \\[8pt]
\Lambda\varphi(\c{L}^-{}^{-a}_{-a})=[\lambda_a, x^a]
&=\left\{ \begin{array}{l}
-\gamma_a \Lambda q k \omega_a^{-1} r_{a-1}^{-1} r_a \\[4pt]
\gamma_{a} \Lambda k \omega_{|a|-1}^{-1} r_{|a|}^{-1} r_{|a|-1}
\end{array}
\right.
&  \begin{array}{l}
\mbox{ for $a>1$,}\\[4pt]
\mbox{ for $a<-1$,}
\end{array} \\[16pt]
\Lambda\varphi(\c{L}^-{}^0_0)=[\lambda_0, x^0]&= -q^{\frac{1}{2}} h \gamma_0
\Lambda &
\mbox{ for $N$ odd} \\[12pt]
\Lambda\varphi(\c{L}^-{}^{-1}_{-1})=[\lambda_1, x^1]&=
\left\{
\begin{array}{l}
-q k \lambda_1 x^1\\[4pt]
-\gamma_1 \Lambda q h (x^0)^{-1} r_1
\end{array}
\right.
& \begin{array}{l}
\mbox{ for $N$ even,}\\[4pt]
\mbox{ for $N$ odd,}
\end{array} \\[16pt]
\Lambda\varphi(\c{L}^-{}^1_1)=[\lambda_{-1}, x^{-1}]&=
\left\{ \begin{array}{l}
-q k \lambda_{-1} x^{-1} \\[4pt]
\gamma_{-1} \Lambda h (r_1)^{-1} x^0
\end{array}
\right.
& \begin{array}{l}
\mbox{ for $N$ even,}\\[4pt]
\mbox{ for $N$ odd,}
\end{array}
\end{array}                                           \label{xlambdarel1}
\ee
To obtain these relations we have used (\ref{xrrel}), (\ref{xLambda})
and (\ref{explicitx}). In the case $a=i$ we also need (\ref{rutil2}),
and for even $N$, $|a|=1$ (\ref{xkapparel}).  It will be noticed that
although complicated in appearence the system of
Equations~(\ref{xlambdarel1}) is actually mainly the first two
equations, which are quite simple, plus a series of special cases when
$i=a$.  The commutation relations between the $x^i$ and the
$\lambda_a$ are independent of the normalization of the latter, so
that they impose no restriction on the $\gamma_a$.

Writing down the explicit expression for the $\hat R$-matrix (\ref{defR})
and using the fact that $[\lambda_a, x^i]=0$ for $i<a$, one finds that
the relation (\ref{Rxlambda}) becomes: 
\be
\begin{array}{lr}
x^h [\lambda_a,x^i]=q[\lambda_a,x^i]x^h+ kq [\lambda_a,x^h]x^i 
&\mbox{ for $h<i$, $h \neq -i$} \nn\\[3mm]
x^h [\lambda_a,x^i]=q[\lambda_a,x^i]x^h &\mbox{ for }
\left\{\begin{array}{l}h>i,h \neq -i,\\h=i=0; \end{array}\right.\\[3mm]
x^i [\lambda_a,x^i]=q^2[\lambda_a,x^i]x^i &\mbox{ for $h=i \neq 0$;}
\end{array}                                         \label{explicitRxlambda1}
\ee
finally, when $h=-i$,
$$
\begin{array}{lr}
x^{-i} [\lambda_a,x^i]=[\lambda_a,x^i] x^{-i} +\\\hskip 2.5cm
kq\Big([\lambda_a,x^{-i}] x^i
- {\displaystyle \sum_{k<i}} q^{-\rho_k+\rho_i}
[\lambda_a,x^k] x^{-k}\Big), &\mbox{ for $i>0$,}\\[6mm]
x^{-i} [\lambda_a,x^i]=[\lambda_a,x^i] x^{-i}
- kq {\displaystyle \sum_{k<i}} q^{-\rho_k+\rho_i} [\lambda_a,x^k]
x^{-k} &\mbox{ for $i<0$.}
\end{array}                                         
$$
These relations can be checked one by one with a
lengthy but straightward calculation by substituting the explicit
expression (\ref{xlambdarel1}) for $[\lambda_a,x^i]$ and commuting
$x^h$ through it.  More particularly, in the case $i<a$ both sides of
the equations are identically 0, because both $[\lambda_a,x^i]$
and $R=P\hat R$ are lower-triangular matrices.  In the case $i>a$ we first
use (\ref{xlambdarel1}) again to commute $x^h$ with $\lambda_a$. Then
we need (\ref{explicitx}) to commute $x^h$ with $x^i$ if $h \neq -i$,
while we need to apply (\ref{rutil1}) to the
expressions of the type $x^k x^{-k}$ to write them in terms of
$r^2_{|k|}$ and $r^2_{|k|-1}$ if $h=-i$. In the case $i=a$ we need
(\ref{xrrel}) to commute $x^h$ through $r_{|a|}$ and $r_{|a|-1}$.  In
the particular case $i=h=0$ (\ref{explicitRxlambda1}) follows from
the $\Lambda x$ commutation relation (\ref{xLambda}).
\ep

Next, look for $\gamma_a$ such that equations (\ref{LLg}),
(\ref{sfilza4}) for $\varphi^-(\c{L}^-)$ with $i=-h$,  
\be
\varphi^-(\c{L}^-{}^i_i\c{L}^-{}^{-i}_{-i})=1, 
\qquad\varphi^-(\c{L}^-{}^0_0)=1,
\ee
are fulfilled. Using (\ref{xlambdarel1}), we easily
find the $\gamma_a$'s given in equations (\ref{gamma}).
\begin{lemma}                                                \label{lulu}
With the $\gamma_a$ given in equations (\ref{gamma}) the 
elements $\lambda_a$ are also solutions to 
equations (\ref{lalarel}) and equations 
\be
\lambda_a [\lambda_b,x^i]=q^{-1} 
(\hat{R}^{-1})^{cd}_{ab} [\lambda_c, x^i]\lambda_d.           \label{Rlambdax}
\ee
\end{lemma}
We shall occasionally use the short-hand notations
\be
e^i_a:=[\lambda_a,x^i]\qquad\qquad\bar e^i_a:=[\bar\lambda_a,x^i]
\label{short-hand}
\ee
Note that, because of (\ref{gRrel}), the $\varphi^-$ images of 
equations (\ref{L-L-rel}) and (\ref{LLg}) 
are respectively equivalent to the `RTT-relations'
\be
\hat R^{ij}_{kl} e_a^k e_b^l= e^i_c e^j_d \hat R^{cd}_{ab}  \label{ree} \\
\ee
and the `gTT-relations'
\be
g^{ab} e^i_a e^j_b= g^{ij} \Lambda^2, \qquad
g_{ij} e^i_a e^j_b= g_{ab} \Lambda^2.                        \label{gtt}
\ee
\begin{prop}                                                \label{prop2}
With the $\gamma_a$ given in equations (\ref{gamma})  
the matrices $e^i_a$ fulfill (\ref{ree}), (\ref{gtt}), 
and the $\varphi^-$ image of (\ref{sfilza4}).
\end{prop}
This will conclude the proof of Theorem \ref{theor1}.

\bp{of Lemma \ref{lulu}}
It is interesting to note that the commutation relations
(\ref{lalarel}) between the $\lambda_a$ are the same as those
(\ref{xrel}) satisfied by the $x^i$, because
$\c{P}_a{}^{ab}_{cd}=\c{P}_a{}^{cd}_{ab}$.  
As (\ref{xrel}) is equivalent to (\ref{explicitx}), so will
equations (\ref{lalarel}) be equivalent to
\be
\begin{array}{ll}
\lambda_a \lambda_b = q \lambda_b \lambda_a & \mbox{ for } a<b, a \neq -b, 
\\[3mm]
[\lambda_a,\lambda_{-a}]= k \omega_{a-1}^{-1} s^2_{a-1} 
&  \mbox{ for } a > 1, \\[4mm]
[\lambda_1, \lambda_{-1}]=\left\{
\begin{array}{l}
0 \\
h s_0^2
\end{array}
\right.&
\begin{array}{l}
\mbox{ for $N$ even,} \\
\mbox{ for $N$ odd,}
\end{array}
\end{array}                                       \label{explicitlambda}
\ee
where the quantities $s^2_a$ are defined by the equation
\be
s^2_a=\sum_{c,d=-a}^a g^{cd} \lambda_c \lambda_d           \label{defs}  
\ee
\pn
for $a \ge 0$ in the case $N$ odd, and for $a\ge 1$ in the
case $N$ even (in the latter case the sum of course runs over $c,d\neq 0$).

It is easy to show the commutation relations (for $i\ge 0$)
\be
r_i \lambda_a=\left\{ 
\begin{array}{lr}
q^2 \lambda_a r_i & \mbox{ for } a<-i, \\[4pt]
q \lambda_a r_i & \mbox{ for } |a| \le i, \\[4pt]
\lambda_a r_i & \mbox{ for } a>i,
\end{array}
\right.                                                     \label{rlambda}
\ee
\be
\lambda_a \Lambda=q^{-1} \Lambda \lambda_a,               \label{lambdaLambda}
\ee
and, for $N$ even, 
\bea
&&[K, \lambda_b]=0 \qquad \mbox{for } |b| \neq 1,     \label{kappalambda} \\
&&K \lambda_{\pm 1}=q^{\mp 1} \lambda_{\pm 1} K.     \nonumber
\eea
which follow from~(\ref{kappaLambda}). 

To show now the relation~(\ref{lalarel})$_1$ we consider first the case
$a<b$, excluding the cases $N$ odd and $a=0$, and $N$ even
and $a=\pm 1$. 
By using (\ref{xlambdarel1})$_2$, (\ref{rlambda}) and 
(\ref{lambdaLambda}) we obtain respectively the identities
\[
\lambda_a \lambda_b =
\gamma_a \Lambda r_{|a|}^{-1} r_{|a|-1}^{-1} \lambda_b x^{-a} =
\gamma_a \Lambda \lambda_b r_{|a|}^{-1} r_{|a|-1}^{-1} x^{-a} = 
q \lambda_b \lambda_a.
\]
If $N$ is even and $a=\pm 1$  we obtain
\be
\lambda_{\pm 1} \lambda_b =
\gamma_{\pm 1} \Lambda (x^{\pm 1})^{-1} \lambda_b  K^{\mp 1}=
\gamma_1 \Lambda \lambda_b (x^{\pm 1})^{-1} K^{\mp 1} = 
q \lambda_b \lambda_{\pm 1},
\ee
using respectively the identities (\ref{kappalambda}),
(\ref{xlambdarel1}) and (\ref{lambdaLambda}). We proceed similarly in
the case $|b|=1$ when $N$ is even.  The calculation is analogous for
the other cases $b \neq -a$.  Summing up, for $b \neq -a$ the
$\lambda_a$ and $\lambda_b$ $q$-commute, so that there is no
restriction on the normalization constants $\gamma_a$.

We now consider the cases $a=-b$.
It follows from (\ref{deflambda}), (\ref{gamma}) that 
\be
\begin{array}{ll}
s^2_0 &=\Lambda^2q^{-2}h^{-2}(x_0)^{-2} \qquad\mbox{ for $N$ odd}\cr
s^2_1 &= \Lambda^2q^{-2}k^{-2}\omega_1^2r_1^{-2}
        \qquad\mbox{ for $N$ even.}           
\end{array}                                            \label{trigger}
\ee
We use these two relations as initial steps to show by induction that
\be
s^2_a=q^{-2}\Lambda^2 \omega_a^2k^{-2}r_a^{-2}
\qquad \mbox{for $a\ge 1$.}                            \label{simpler}
\ee
In fact
\bea
s^2_a &\stackrel{(\ref{defs})}{=}& s^2_{a-1}+q^{-\rho_a}\lambda_a
\lambda_{-a}+q^{\rho_a}\lambda_{-a}\lambda_a \nn
&\stackrel{(\ref{deflambda})}{=}& s^2_{a-1}+\Lambda^2\gamma_a
\gamma_{-a}r_a^{-2}r_{a-1}^{-2}\left[q^{-2-\rho_a}x^{-a}x^a+
q^{\rho_a}x^ax^{-a}\right]  \nn
&\stackrel{(\ref{rutil1})}{=}& s^2_{a-1}+\Lambda^2\gamma_a
\gamma_{-a}r_a^{-2}r_{a-1}^{-2}\left[q^{-2-\rho_a}(\omega_a^{-1}r_a^2-
q\omega_{a-1}^{-1}r_{a-1}^2)\right.\nn
&& \left.+ q^{\rho_a}(\omega_a^{-1}r_a^2-
q^{-1}\omega_{a-1}^{-1}r_{a-1}^2)\right] \nn
&=& s^2_{a-1}+\Lambda^2\gamma_a
\gamma_{-a}r_a^{-2}r_{a-1}^{-2}\left[q^{-1}\omega_{a-1}\omega_a^{-1}r_a^2
-q^{-1}\omega_a\omega_{a-1}^{-1}r_{a-1}^2\right] \nn
&\stackrel{(\ref{gamma})}{=} &s^2_{a-1}-\Lambda^2q^{-2}k^{-2}
\omega_{a-1}^2r_{a-1}^{-2}+\Lambda^2q^{-2}k^{-2}
\omega_a^2r_a^{-2}. \nonumber
\eea
Assuming that (\ref{simpler}) holds for $a=b-1$, the first two terms
in the last line are opposite and therefore cancel, and the third
gives (\ref{simpler}) for $a=b$, as claimed.  We now consider the
commutators $[\lambda_a,\lambda_{-a}]$ with $a\ge 1$. For $N$ even we
find the claim $[\lambda_1,\lambda_{-1}]=0$ by a straightforward
calculation.  In all other cases we proceed as follows,
\bea
[\lambda_a,\lambda_{-a}]&\stackrel{(\ref{deflambda}),(\ref{xrrel})}{=}&
\gamma_a \gamma_{-a}q^{-1}\Lambda^2r_a^{-2}r_{a-1}^{-2}(q^{-1}x^{-a}
x^a-qx^ax^{-a}) \nn
&\stackrel{(\ref{rutil2})}{=}&-q^{-1}\Lambda^2\gamma_a \gamma_{-a}
\omega_a^{-1}kr_{a-1}^{-2}
\stackrel{(\ref{gamma})}{=}q^{-2}k^{-1}\omega_{a-1}
\Lambda^2r_{a-1}^{-2}\nn
&\stackrel{(\ref{simpler})}{=}&k \omega_{a-1}^{-1} s^2_{a-1}, \nonumber
\eea
as claimed.

With the choice (\ref{gamma}) for the normalization constants $\gamma_a$,
the algebra generated by $x^i, \lambda_i, \Lambda^{\pm 1}, K^{\pm 1},
r_i^{\pm 1}$ is symmetric, i.e. is invariant, with respect to the
following transformation $\c{S}$ 
\be
\begin{array}{lllr}
\lambda_{\pm 1}& \rightarrow &
\left( \gamma_{\pm 1}(q^{-1})\gamma_{\mp 1}(q)^{-1} \right)^{\frac{1}{2}}
x^{\mp 1} & \mbox{for $N$ even}, \\ \\
\lambda_a& \rightarrow & q^{-\frac{1}{2}}
\left( \gamma_a(q^{-1}) \gamma_{-a}(q)^{-1} \right)^{\frac{1}{2}} x^{-a} &
\mbox{otherwise,} 
 \\ \\
x^{\pm 1}& \rightarrow &
\left( \gamma_{\pm 1}(q^{-1})\gamma_{\mp 1}(q)^{-1}\right)^{\frac{1}{2}}
\lambda_{\mp 1}  & \mbox{for $N$ even}, \\ \\
x^a& \rightarrow &q^{-\frac{1}{2}}
\left( \gamma_a(q^{-1})\gamma_{-a}(q)^{-1} \right)^{\frac{1}{2}} \lambda_{-a} 
& \mbox{otherwise,} \\ \\
\Lambda& \rightarrow &\Lambda, \\ \\
K& \rightarrow &K^{-1}, \\ \\
q & \rightarrow & q^{-1}.
\end{array}
\label{symmetry}
\ee
Notice that $\c{S}$ is an involution, $\c{S}^2=\mbox{id}$.
Because of (\ref{propR1}), (\ref{propR2}) and (\ref{propR3}) under
$\c{S}$ the $xx$-commutation relations (\ref{xrel}) and the
$\lambda \lambda$-commutation relations (\ref{lalarel}) are exchanged.
The $x^a \Lambda$ (\ref{xLambda}) and the $\lambda_a \Lambda$ relations
(\ref{lambdaLambda}) are exchanged as well, while the $x \lambda$ relations
(\ref{xlambdarel1}) are invariant. We can immediately check that for $N$ odd
under $\c{S}$
\be
r_a^2 \rightarrow \sum_{b=-a}^{a} q^{\rho_b} q^{-1} \lambda_{-b} \lambda_{b}
(\gamma_b(q^{-1}) \gamma_{-b}(q^{-1}) \gamma_b(q)^{-1}
\gamma_{-b}(q)^{-1})^{\frac{1}{2}}=s_a^2 
\ee
where we have used (\ref{gamma}). 
In the case of even $N$ the same calculation holds for the terms with $|b|>1$,
but we have to treat the term with $|b|=1$ separately
\be
r_1^2 \rightarrow 2 \lambda_{-1} \lambda_{1}
(\gamma_1(q^{-1}) \gamma_{-1}(q^{-1}) \gamma_1(q)^{-1}
\gamma_{-1}(q)^{-1})^{\frac{1}{2}}=s_1^2 
\label{s1}
\ee
The relation (\ref{deflambda}) expressing $\lambda_a$ in terms of $x^{-a}$
is invariant as well. To see this we first 
express $r_a^{-1}$ and $r_{a-1}^{-1}$ in (\ref{deflambda}) through
$s_a$ and $s_{a-1}$ respectively, then we use
(\ref{lambdaLambda}) to move $\Lambda^{-1}$ to the left.
In this way we are able to rewrite (\ref{deflambda}) in the form
\be
\lambda_a=\gamma_a(q) q^2 \Lambda^{-1} s_a s_{a-1} x^{-a}.
\label{lambdabis}
\ee
Taking into account (\ref{gamma}) and (\ref{xLambda}), under $\c{S}$
(\ref{lambdabis}) becomes
\bea
x^{-a} &=& \gamma_a(q)^{-1} r_a r_{a-1} \Lambda^{-1} \lambda_a
\eea
i.e. we recover (\ref{deflambda}).
Again, in the case of even $N$ the special case $|a|=1$ has to be
treated separately, but due to (\ref{symmetry})${}_6$ and (\ref{s1}),
it is easily checked that (\ref{deflambda}) is
invariant in this case, too.

This transformation is useful, because it enables us to get (\ref{Rlambdax})
by applying $\c{S}$ to (\ref{Rxlambda}) and then using the
properties (\ref{propR1}), (\ref{propR2}) and (\ref{propR3}) of
the $\hat R$-matrix and (\ref{gamma}). For $N$ odd and
$N$ even, $h \neq -i$
\bea
\lefteqn{\lambda_h [x^a,\lambda_i]} \nn
&=&\sum\limits_{j,k}q^{-1}
\hat R(q^{-1}){}^{-h-i}_{-j-k}
\sqrt{\frac{\gamma_h(q)}{\gamma_{-h}(q^{-1})}
\frac{\gamma_i(q)}{\gamma_{-i}(q^{-1})} 
\frac{\gamma_{-j}(q^{-1})}{\gamma_j(q)}
\frac{\gamma_{-k}(q^{-1})}{\gamma_k(q)}} [x^a,\lambda_j]
\lambda_k \nn
&=&\sum\limits_{j,k}q^{-1} \hat R(q)^{-1} {}^{jk}_{hi} [x^a,\lambda_j] 
\lambda_k.  
\eea
In the particular case that $N$ is even and $h=-i$ from
\begin{eqnarray*}
\gamma_1(q) \gamma_{-1}(q)&=&\gamma_{1}(q^{-1})\gamma_{-1}(q^{-1}) \\
\gamma_b(q) \gamma_{-b}(q)&=&q^{-2}\gamma_b(q^{-1}) \gamma_{-b}(q^{-1})
\quad \mbox{ for } b \neq \pm 1
\end{eqnarray*}
and the property (\ref{propR3}) of the $\hat R$-matrix,
it is easily seen that (\ref{Rlambdax}) still holds. This concludes
the proof of Lemma \ref{lulu}.
\ep

\bp{of Proposition \ref{prop2}}
To prove (\ref{ree}), we use (\ref{Rlambdax}) and (\ref{Rxlambda})
\bea
&& \hat R^{cd}_{ab} [\lambda_c, x^i][\lambda_d, x^j]=
 \hat R^{cd}_{ab} (\lambda_c x^i-x^i \lambda_c)[\lambda_d, x^j]= \nn
&& \hat R^{cd}_{ab} (q \hat R^{ij}_{kl} \lambda_c
[\lambda_d, x^k] x^l - q^{-1}(\hat R^{-1})^{ef}_{cd} x^i
[\lambda_e, x^j] \lambda_f)=\nn
&&
\hat R^{cd}_{ab} q^{-1}
(\hat R^{-1})^{ef}_{cd} q \hat R^{ij}_{kl} [\lambda_e,x^k]
(\lambda_f x^l - x^l \lambda_f)=
\hat R^{ij}_{kl} [\lambda_a, x^k] [\lambda_b,x^l],
\nonumber
\eea
i.e. the `RTT'-relations for $e^i_a$.
By repeated application of the `RTT'-relations it is an immediate
result that for any polynomial $f(\hat R)$ 
\be
f(\hat R^{ij}_{kl}) e^k_a e^l_b= e^i_c e^j_df (\hat R^{cd}_{ab}).
\ee
In particular the projectors $\c{P}_s$, $\c{P}_a$, $\c{P}_t$ are of this
form. If we write $\c{P}_t$ explicitly using (\ref{Pt}), this yields the
$gTT$-relations (\ref{gtt}) 
[which are equivalent to the $\varphi^-$-image of equations 
(\ref{LLg})] also for $h\neq -i$, which we had not proved yet.
\ep

\subsection{Proof of Theorem \ref{theor2}}

The proof of \ref{theor2} is similar to the one of Theorem 
\ref{theor1}. The explicit expressions for
$\Lambda^{-1} \varphi^+(\c{L}^+{}^{-a}_{-i})$ are
\be
\begin{array}{lll}
\Lambda^{-1} \varphi(\c{L}^+{}^{-a}_{-i})
=q^{\rho_a-\rho_i} [\bar \lambda_a, x^i] &=0& \mbox{ for $i>a$,} \\[4pt]
\Lambda^{-1}\varphi(\c{L}^+{}^{-a}_{-i})
=q^{\rho_a-\rho_i} [\bar \lambda_a, x^i]&=q^{\rho_a-\rho_i-1} k
\bar \lambda_a x^i& \mbox{ for $i<a$,} 
\\[8pt]
\Lambda^{-1} \varphi(\c{L}^+{}^{-a}_{-a})=[\bar \lambda_a, x^a]
&=\left\{ \begin{array}{l}
-\bar \gamma_a \Lambda^{-1} k \omega_{a-1}^{-1} r_a^{-1} r_{a-1}
\\[4pt]
\bar \gamma_a \Lambda^{-1} k q^{-1} \omega_a^{-1} r_{|a|-1}^{-1} r_{|a|}
\end{array}
\right.
&  \begin{array}{l}
\mbox{ for $a>1$,}\\[4pt]
\mbox{ for $a<-1$,}
\end{array} \\[16pt]
\Lambda^{-1} \varphi(\c{L}^+{}^0_0)=[\bar \lambda_0, x^0]&= 
\bar \gamma_0 \Lambda^{-1} q^{-\frac{1}{2}} h 
& \mbox{ for $N$ odd} \\[12pt]
\Lambda^{-1}\varphi(\c{L}^+{}^{-1}_{-1})=[\bar \lambda_1, x^1]&=
\left\{
\begin{array}{l}
q^{-1} k \bar \lambda_1 x^1\\[4pt]
-\bar \gamma_1 \Lambda^{-1} h r_1^{-1} x^0
\end{array}
\right.
& \begin{array}{l}
\mbox{ for $N$ even,}\\[4pt]
\mbox{ for $N$ odd,}
\end{array} \\[16pt]
\Lambda^{-1}\varphi(\c{L}^+{}^1_1)=[\bar \lambda_{-1}, x^{-1}]&=
\left\{ \begin{array}{l}
q^{-1} k \bar \lambda_{-1} x^{-1} \\[4pt]
\bar \gamma_{-1} \Lambda^{-1} h q^{-1} r_1 (x^0)^{-1}
\end{array}
\right.
& \begin{array}{l}
\mbox{ for $N$ even,}\\[4pt]
\mbox{ for $N$ odd.}
\end{array}
\end{array}                                           \label{bareia}
\ee

\subsection{Proof of Theorem \ref{theor3}}            \label{proofs2}

Using relations (\ref{gRrel}), it is easy to show that
the image under $\varphi$ of (\ref{L+L-rel}) is equivalent to
\be
\hat R^{cd}_{ab}\bar e_c^i e_d^j= \hat R^{ij}_{kl}e_a^k\bar e_b^l
\label{equiv}
\ee
in the notation (\ref{short-hand}).

Theorem \ref{theor1} (\ref{theor2}) fixes the coefficients
$\gamma_a$ ($\bar \gamma_a$) for $a<0$ in terms of 
$\gamma_a$ ($\bar \gamma_a$) for $a \ge 0$.
We can use the remaining freedom in the choice of $\gamma_a$, $\bar \gamma_a$
to find further conditions on $\gamma_a$ for $a>0$,
and relations relating the coefficients $\bar \gamma_a$ to $\gamma_a$  
so that (\ref{sfilza3}) and (\ref{L+L-rel})
are fulfilled. We start with the observation that 
in the case of odd $N$ (\ref{deflambda}) and (\ref{defbarlambda}) imply
\be
\bar \lambda_a=\bar \gamma_a \gamma_a^{-1} \Lambda^{-2} \lambda_a 
\label{blambda}
\ee
and use the equations (\ref{Rlambdax}) and (\ref{Rxlambda}). 
In this way we see that
\bea
&&\sum\limits_{c,d}\hat R^{cd}_{ab} [\bar \lambda_c, x^i][\lambda_d, x^j]=
\sum\limits_{c,d}\hat R^{cd}_{ab} \bar \gamma_c \gamma_c^{-1} 
(\Lambda^{-2} \lambda_c x^i-x^i \Lambda^{-2} \lambda_c) [\lambda_d, x^j]= \nn
&&\sum\limits_{c,d}\hat R^{cd}_{ab} \bar \gamma_c \gamma_c^{-1} \Lambda^{-2} 
(\lambda_c x^i -q^{-2} x^i \lambda_c) [\lambda_d, x^j]= \nn 
&&\sum\limits_{c,d}\hat R^{cd}_{ab} \bar \gamma_c \gamma_c^{-1} \Lambda^{-2}
(q \sum\limits_{k,l}\hat R^{ij}_{kl} \lambda_c [\lambda_d, x^k] x^l
- q^{-3} \sum\limits_{e,f}\hat R^{-1}{}^{ef}_{cd} x^i [\lambda_e, x^j] 
\lambda_f)=\nn
&&\sum\limits_{c,d}\sum\limits_{e,f}\sum\limits_{k,l}\hat R^{cd}_{ab} 
\hat R^{-1}{}^{ef}_{cd} \hat R^{ij}_{kl} \bar \gamma_c \gamma_c^{-1} 
[\lambda_e,x^k](\Lambda^{-2}\lambda_f x^l - x^l \Lambda^{-2} \lambda_f)=\nn
&& \sum\limits_{c,d}\sum\limits_{e,f}\sum\limits_{k,l}\hat R^{cd}_{ab} 
\hat R^{-1}{}^{ef}_{cd} \hat R^{ij}_{kl} \bar \gamma_c \gamma_c^{-1} \gamma_f 
\bar \gamma_f^{-1}[\lambda_e, x^k] [\bar \lambda_f,x^l], \nonumber
\eea
If the further condition holds,
\be
\bar \gamma_a \gamma_a^{-1}=c\equiv\mbox{const,}     \label{cost}
\ee
then in the last line the $\gamma$'s cancel with each other,
so do $\hat R$ and $\hat R^{-1}$, and we find that (\ref{equiv}) is
actually satisfied.

To see that (\ref{cost}) is not only a sufficient, but also a necessary
condition for (\ref{equiv}),
we write down the latter for the particular values of the indices
$i=a$, $j=b=a+1$ for $a=-n, \ldots n-1$.
\bea
[\bar \lambda_{a+1}, x^a] [\lambda_a, x^{a+1}]+k[\bar \lambda_a,
x^a][\lambda_{a+1},x^{a+1}]=\label{equivaa}\\
{[\lambda_a,x^{a+1}]} [\bar \lambda_{a+1}, x^a]+k[\lambda_a, x^a][\bar
\lambda_{a+1},x^{a+1}]
\nonumber
\eea
We plug in the expressions
(\ref{xlambdarel1}),(\ref{bareia}) for $e^i_a$ and $\bar e^i_a$
and apply the relations (\ref{rutil1}). For $a>1$, (\ref{equivaa})
implies
\be
k^3 q (\bar \gamma_a \gamma_{a+1}- \bar \gamma_{a+1} \gamma_a) r_{a-1}
r_{a+1}r_a^{-2}\omega_{a-1}^{-1}\omega_{a+1}^{-1}=0
\ee
This means that in order for (\ref{equiv}) to hold, we have to require
\be
\bar \gamma_a \gamma_a^{-1}= \bar \gamma_{a+1} \gamma_{a+1}^{-1}
\label{nec}
\ee
for every $a>1$. A similar reasoning can be repeated for $a \le 1$
to show that (\ref{nec}) has to hold for any value of $a$.
Therefore (\ref{cost}) is necessary.

When $N$ is even, it is not possible to satisfy (\ref{equiv}).
This is a consequence of the fact that (\ref{blambda}) does not hold for
$|a|=1$ in this case, due to the particular form of $\lambda_{\pm 1}$ and
$\bar \lambda_{\pm 1}$. Choose the indices in (\ref{equiv}) to be e.g.
$i=a=1$, $j=b=2$. Plug in (\ref{xlambdarel1}), (\ref{bareia}) for $e^i_a$
and $\bar e^i_a$ and the definitions (\ref{deflambda}), (\ref{defbarlambda})
for $\lambda_a$ and $\bar\lambda_a$, then apply (\ref{xkapparel}) and
(\ref{explicitx}). In this way (\ref{equivaa}) becomes
\begin{eqnarray*}
-k^2 \left( \bar\gamma_2 \gamma_1 r_2^{-1} r_1^{-1} x^{-2} x^2 K^{-1}+
k \bar \gamma_1 \gamma_2 q \omega_2^{-1} r_2 r_1^{-1} K\right)=\\
-k^2 \bar\gamma_2 \gamma_1 r_2^{-1} r_1^{-1}
\left( x^2 x^{-2}-k \omega_1^{-1} r_1^2 \right) K^{-1}
\end{eqnarray*}
Due to the commutation relation (\ref{explicitx}) between $x^2$ and $x^{-2}$
the terms which are proportional to $K^{-1}$ cancel and the following
equation should be satisfied
$$
k^3 q \bar \gamma_1 \gamma_2 \omega_2^{-1} r_2 r_1^{-1} K=0.
$$
But this would mean that either $\bar \gamma_1$ or $\gamma_2$ should
vanish, which is not possible, if we want to have $N$ independent objects
$\lambda_a$ instead of fewer. That is the reason why the theorem does
not hold for $N$ even.

In the case $N$ odd (\ref{cost}) is consistent with (\ref{gamma}), 
(\ref{bargamma}). We can determine $c$ e.g. by applying 
(\ref{cost}) to $a=0$ and by recalling (\ref{gamma})$_1$, 
(\ref{bargamma})$_1$. We thus find $c=-q$.
In this way we get the last of the equations (\ref{fixgamma}),
which completely fixes the coefficients $\bar \gamma_a$ in terms of
the $\gamma_a$. 

As the last step, we require (\ref{sfilza3}). This imposes another condition
relating $\bar \gamma_a$ to $\gamma_a$:
\bea
&\bar \gamma_a \gamma_a= k^{-2} \omega_{a-1} \omega_a q^{-1} &
\mbox{for } a>1,                                         \label{sfilza3a1} \\
&\bar \gamma_1 \gamma_1= h^{-2} q^{-1} \label{sfilza3a2}\\
& \bar \gamma_0 \gamma_0= -h^{-2} \label{sfilza3a3}\\
&\bar \gamma_{-1} \gamma_{-1}= h^{-2} q \label{sfilza3a4}\\
& \bar \gamma_a \gamma_a= k^{-2} \omega_{a-1} \omega_a q &
\mbox{for } a<-1, \label{sfilza3a5}
\eea
Here we have used the expression (\ref{xlambdarel1}) for 
$e^a_a=\Lambda\varphi(\c{L}^-{}^{-a}_{-a})$ and (\ref{bareia})
for $\bar e^a_a=\Lambda\varphi(\c{L}^+{}^{-a}_{-a})$
Equations (\ref{sfilza3a1}) to (\ref{sfilza3a5}) are compatible with
(\ref{gamma}), (\ref{bargamma}).  If we replace 
$\bar\gamma_a=-q\gamma_a$ we find the remaining equations in 
(\ref{fixgamma}).


\providecommand{\href}[2]{#2}\begingroup\raggedright\endgroup

\end{document}